\documentclass[11pt,a4paper,twoside]{article}
\usepackage{amsmath, exscale, amssymb, graphicx}
\small
\normalsize
\usepackage[latin1]{inputenc}
\usepackage[T1]{fontenc}
\usepackage{amsfonts}
\usepackage{amsbsy}
\usepackage{fancyheadings}
\usepackage{amscd}
\usepackage{theorem}
\usepackage{multicol}
\allowdisplaybreaks

\def\build#1_#2^#3{\mathrel{\mathop{\kern 0pt#1}\limits_{#2}^{#3}}}
\def\noi{{\noindent}}
\def\be{\begin{equation}}
\def\ee{\end{equation}}
\def\ba{\begin{eqnarray*}}
\def\ea{\end{eqnarray*}}

\def\E{{\bf E}}
\def\cqfd{ \hfill $\blacksquare$ }

\def\llb{[\hspace{-.10em} [ }
\def\rrb{ ] \hspace{-.10em}]}

\def\g{{\cal G}}
\def\aal{{\bf \alpha}}

\def\f{{\cal F}}

\def\m{{\cal M}}
\def\h{{\cal H}}
\def\ii{{\cal I}}
\def\I{{\cal I}}

\def\t{{\cal T}}
\def\v{{\cal V}}

\def\P{{\bf P}}
\def\EE{{\bf E}}%\def\MM{{\bf M}}

\def\ee{\epsilon}

\def\dist{{\bf d}}

\def\noi{\noindent}

\def\build#1_#2^#3{\mathrel{\mathop{\kern 0pt#1}\limits_{#2}^{#3}}}

\newtheorem{theorem}{Theorem}[section]
\newtheorem{lemma}[theorem]{Lemma}
\newtheorem{proposition}[theorem]{Proposition}

{\theorembodyfont{\rmfamily}\newtheorem{remark}{Remark}[section]}

{\theorembodyfont{\rmfamily}\newtheorem{example}{Example}[section]}
\newcommand{\R}{\mathbb{R}}
\newcommand{\RR}{\mathbb{R}}

\newcommand{\Z}{\mathbb{Z}}

\newcommand{\Q}{\mathbb{Q}}

\newcommand{\D}{\mathbb{D}}
\newcommand{\N}{\mathbb{N}}
\newcommand{\U}{\mathbb{U}}
\newcommand{\F}{\mathbb{F}}

\newcommand{\M}{\mathbb{M}}

\newcommand{\T}{\mathbb{T}}
\newcommand{\un}{\boldsymbol{1}}

\begin{document}

\title{ {\bf   CONTINUUM RANDOM TREES  \\
AND BRANCHING PROCESSES WITH IMMIGRATION} }
\author{ by \\
Thomas {\sc Duquesne}
 \thanks{  Université Paris-Sud, Mathématiques, 91405 Orsay Cedex,
   France; email: thomas.duquesne@math.u-psud.fr; 
supported by NSF Grants DMS-0203066 
and DMS-0405779} }
\vspace{4mm}
\date{\today} 

\maketitle

\begin{abstract}

We study a genealogical model for continuous-state branching processes with
immigration with a (sub)critical branching
mechanism. This model allows  
the immigrants to be on the same line 
of descent. The corresponding family tree is an ordered rooted 
continuum random tree
with a single infinite end defined by two continuous processes denoted
by $(\overleftarrow{H}_t ;t\geq 0 )$ and 
$(\overrightarrow{H}_t ;t\geq 0 )$ that code the parts at resp. the left and
the right hand side 
of the infinite line of descent of the tree. These processes are called the left and the right height
processes. 
We define their local time processes via
an approximation procedure 
and we prove that they enjoy a Ray-Knight property.
We also discuss the important 
special case corresponding to 
the size-biased Galton-Watson tree in the continuous setting.
%and we show that, as in the discrete setting, 
%they are the limit of the 
%genealogy of (sub)critical L\'evy trees conditioned on non-extinction. 
In the last part of the paper we give a convergence result under general 
assumptions for rescaled discrete 
left and right contour processes of sequences of Galton-Watson trees with immigration.
We also provide a strong invariance principle for a sequence of 
rescaled Galton-Watson processes with immigration that also holds in the
supercritical case.

\vspace{4mm}

\noi {\it MSC 2000 subject classifications:} 60F17,60G17, 60J80

\vspace{4mm}

\noi {\it Key words and phrases:}
continuous-state branching process, immigration,
size-biased tree, L\'evy tree, continuum random tree, height process, contour
process, Ray-Knight theorem, limit theorem, invariance principle,
Galton-Watson branching process.

\end{abstract}

\section{Introduction}

\subsection{The genealogy of Galton-Watson branching processes.}
The continuous-state branching processes with immigration (CSBPI for short)
have been introduced by Kawasu and Watanabe in \cite{KaWa}. 
They are the continuous analogues of the Galton-Watson processes with 
immigration. In this paper we discuss a genealogical model for 
CSBPI's that can be described in the discrete setting as follows: let $\mu$ and $\nu$ be two
probability measures on the set of non-negative integers denoted by $\N$. 
Recall that a  Galton-Watson process $Z=(Z_n; n\geq 0)$ 
with {\it offspring distribution} $\mu $ and {\it immigration distribution} 
$\nu $ (a GWI$(\mu, \nu)$-process for short)
is an $\N$-valued Markov chain 
whose transition probabilities are characterized by 
\begin{equation}
\label{transGWI}
\E \left[ x^{Z_{n+1} } \mid Z_n \right]
= g (x)^{Z_{n} } f\left( x\right)\; , \quad x\in [0,1], 
\end{equation} 
where $g$ (resp. $f$) stands for the generating function of $\mu $
(resp. $\nu$).

The genealogical model we consider can be informally described as follows.
Consider a population evolving at random 
roughly speaking according to a GWI$(\mu, \nu)$-process:
The population can be decomposed in two kinds of individuals, namely
the mutants and the non-mutants; there is exactly one mutant at a given generation; 
each individual gives
birth to an independent number of children: the mutants 
in accordance with $\nu$ and the non-mutants in accordance with  $\mu$.  
We require that all the mutants are on the same infinite line of descent.
Except in one part of the paper we restrict our attention to a critical
or subcritical offspring distribution $\mu$:

$$ \bar{\mu}:= \sum_{k\geq 0} k\mu (k) \leq 1 .$$

Then the resulting family tree is a tree with a single infinite
end. We call such trees {\it sin-trees} following Aldous' terminology in 
\cite{Alfringe} in Section 4. In order to code them by real-valued
functions, all the discrete trees that we consider are
ordered and rooted or equivalently are 
planar graphs (see Section 2.1 for more details). So 
we have to specify for each mutant how many of
its children are on the left hand of the infinite line of descent and how
many are on the right hand. We choose to dispatch them independently at
random in accordance with a probability measure $r$ on the set 
$\{ (k,j)\in \N^* \times \N^* \; :\; 1\leq j\leq k \}$. More
precisely, with probability $r(k,j)$  a mutant has 
$k-1$ non-mutant children and one mutant child who is the $j$-th child
to be born; consequently there is $j-1$ non-mutants children on the
left hand of the infinite line of descent and $k-j$ on the right
hand. The immigration distribution $\nu $ then is given by 
$$ \nu (k-1)= \sum_{1 \leq j\leq k} r(k,j) \; ,\quad k\geq 1 .$$
The probability measure $r$ is called 
the {\it dispatching distribution} and the resulting random tree is
called a ($\mu,r$)-{\it Galton-Watson tree with immigration} (a 
GWI$(\mu, r)$-tree for short). Indeed, 
if we denote by $Z_n$ the number of non-mutants 
at generation $n$ in the tree, then it is easy to see that $Z=(Z_n; n\geq 0)$ is 
GWI$(\mu, \nu)$-process.
%(see Figure 1 and 
%Section 2.1 for a precise definition).

Let us mention the special case $r(k,j)=\mu (k)/\bar{\mu} $, $1\leq j\leq k$ 
that corresponds to the size-biased Galton-Watson tree with offspring
distribution $\mu $. This random tree arises naturally by
conditioning  (sub)critical GW-trees 
on non-extinction: see \cite{AlPit98}, \cite{Alfringe}, \cite{Kes87}, 
\cite{Lamb1}, 
\cite{LyoPemPer} and in the continuous case \cite{Li1}.

%\begin{figure}[ht]
%\psfrag{GW}{{\tiny GW}}
%\psfrag{GW1}{{\tiny GW(1)}}
%\psfrag{GW2}{{\tiny GW(2)}}
%\psfrag{GWj}{{\tiny GW(j-1)}}
%\psfrag{GWjj}{{\tiny GW(j+1)}}
%\psfrag{GWk}{{\tiny GW(k)}}
%\psfrag{TAU}{{\tiny infinite line of descent}}
%\psfrag{ROOT}{ {\tiny root} }
%\epsfxsize=5cm
%\centerline{\epsfbox{GWI.eps}}
%\caption{{\small The GWI$(\mu , r )$-tree $\tau$. The bold line represents the
%    infinite line of descent.}}
%\label{GWItree.fig}
%\end{figure}

We shall code a GWI$(\mu, r)$-tree $\tau$ 
by two real-valued functions in the
following way: think of $\tau$ as a planar
graph embedded in the clockwise oriented half-plane with unit edge length and consider a particle visiting 
continuously the edges of 
$\tau$ 
at speed one from the left to the right, going backward as less as possible;
we denote by  $ \overleftarrow{C}_s (\tau ) $ 
the distance from the root of the particle at time $s$
and we call the resulting process $\overleftarrow{C}(\tau):=
(\overleftarrow{C}_s (\tau) ; s\geq 0 )$ 
the {\it left contour 
process} of $\tau$. It is clear that the particle never reaches the right part
of $\tau$ but $\overleftarrow{C} (\tau)$ completely codes the left part of
$\tau$. 
We denote by 
$\overrightarrow{C}(\tau)$ 
the process corresponding to a particle visiting 
$\tau$ from the right to the left so we can then reconstruct $\tau$
from $(\overleftarrow{C}(\tau), \overrightarrow{C} (\tau))$: see Section 2.2 for precise definitions and other codings of sin-trees.

\subsection{Background on continuous-states branching processes with immigration.}
The main purpose of the paper is to
provide a genealogical model for CSBPI's and to build a continuous
family tree coded by two functions playing the  role 
of $\overleftarrow{C}$ and $\overrightarrow{C}$.
Before discussing it more specifically, let us recall from 
\cite{KaWa} that continuous-state
branching processes with immigration are $[0, \infty]$-valued stochastically
continuous Markov processes whose distribution is characterized by two
functions 
on $[0, \infty)$: a {\it branching mechanism}
$\psi$ such that $(-\psi)$ is the Laplace exponent of 
a spectrally positive L\'evy process denoted by 
$X=(X_t; t\geq 0)$ and an {\it immigration mechanism} $\varphi$ that is 
the Laplace exponent of a subordinator denoted by $W=(W_t; t\geq 0)$:
$$ \E \left[ e^{-\lambda X_t}\right]= e^{t\psi (\lambda )} \; , \quad 
\E \left[ e^{-\lambda W_t}\right]= e^{-t\varphi (\lambda )} 
\; , \quad \lambda , t \geq 0 . $$
More precisely, $Y^*=(Y^*_t; t\geq 0)$ is a 
$(\psi,\varphi)$-continuous-state branching process with immigration 
(a CSBPI$(\psi,\varphi)$ for short) if its 
transition kernels are characterized by
\begin{equation}
\label{CSBPImarkov}
\E \left[  e^{-\lambda Y^*_{a+b}} \left|Y^*_b \right.\right]=
 \exp \left( - u (a,\lambda ) Y^*_{b} -\int_0^a ds \, 
\varphi ( u (s,\lambda )) \right) \; , \; a,b,\lambda\geq 0,
\end{equation}
where $u (a,\lambda )$ is the unique nonnegative 
solution of the differential equation
\begin{equation}
\label{equavv}
\frac{\partial}{\partial a}  u (a,\lambda ) = -\psi 
\left( u (a,\lambda ) \right)
\quad {\rm and } \quad  u (0, \lambda ) =\lambda \; ,\; a, \lambda\geq 0.
\end{equation}
Note that this differential equation is equivalent to the integral equation
\begin{equation}
\label{eqCSBP}
\int_{ u (a, \lambda )}^{\lambda } \frac{du}{\psi (u)} = a .
\end{equation}
Observe that $\infty $ is an absorbing state. In the paper, we will only consider 
{\it conservative} processes, that is processes such that a.s. $Y^*_t < \infty$, $t\geq 0$. This is 
equivalent to the analytical conditions
\begin{equation}
\label{conservative}
\int_{0+} \frac{ du}{|\psi (u) |} = \infty 
\quad {\rm and } \quad \psi(0)=\varphi (0)=0 \; . 
\end{equation}
Observe that if there is no immigration, that is if $\varphi =0$, the 
process is simply a 
$\psi$-continuous-state branching process (a CSBP($\psi$) for short). 
We shall 
denote CSBPI's and 
CSBP's in a generic way by resp. $Y^*=(Y^*_t; t\geq 0)$ and 
$Y=(Y_t; t\geq 0)$. We refer to \cite{KaWa} or  \cite{Pi} for a precise
discussion of CSBPI's and 
to \cite{Bi2} and \cite{Grey} for results on CSBP's.

Except in Section 4.1 we only consider CSBPI's with (sub)critical branching
mechanism. This assumption is equivalent to the fact that $X$ does not drift 
to $+\infty $. In that case $\psi $ is of the form 
\begin{equation}
\label{subcri}
 \psi (\lambda )= \alpha \lambda + \beta \lambda^2 +
\int_{(0, \infty)} \pi (dr) (e^{-\lambda r} -1+\lambda r) \; , \quad 
\lambda\geq 0,  
\end{equation}
where $\alpha, \beta \geq 0$ and $\pi $ is a $\sigma $-finite measure 
on $(0, \infty)$ such that
$\int_{(0, \infty)} \pi (dr)(r\wedge r^2) < \infty $. 
We also assume 
\begin{equation}
\label{Hcont}
\int_{ 1}^{\infty } \frac{du}{\psi (u)} < \infty ,
\end{equation}
which is equivalent to the a.s. 
extinction of CSBP($\psi$) (see \cite{Grey} for details). Note that (\ref{Hcont}) implies 
at least one
of the two conditions: $\beta>0$ or $\int_0^1 r\,\pi(dr)=\infty$ that guarantee
that $X$ has infinite variation sample paths (see \cite{Be} for details).

   We build the family tree corresponding to a CSBPI$(\psi,\varphi)$ 
thanks to two continuous processes 
$(\overleftarrow{H}_t ;t\geq 0 )$ and 
$(\overrightarrow{H}_t ;t\geq 0 )$ called  the {\it left and the right height
processes} that are viewed as the 
contour processes of the parts at the left and at the right hand 
of the infinite line of descent.
More precisely, our construction relies on two auxilliary processes:

\noindent
$\bullet $ The first one is the {\it height process} $H=(H_t; t\geq 0)$
introduced by Le Gall and Le Jan \cite{LGLJ1} coding the genealogy of
$\psi$-continuous-state branching processes (see also \cite{DuLG} for
related results). $H$ is obtained as a
functional $H(X)$ of the spectrally 
positive L\'evy process $X$ with exponent $\psi$.  More precisely, 
for every $t\geq 0$, $H_t$ ``measures'' in a local time sense 
the size of the set $\{s\leq t:X_{s-}=\inf_{[s,t]}X_r\}$ (see Section
3.1 for a more precise definition).
Assumption (\ref{Hcont}) is equivalent for 
$H$ to have a continuous
modification. {\it From now on,  
we only consider this modification}.

  An important role is played by the 
excursion measure $N$ of $X$
above its minimum process. 
In the quadratic branching case $\psi(u)=c\,u^2$, $X$  is a (scaled) 
Brownian motion, the height process
$H$ is
a reflected Brownian motion and the ``law'' of $H$ under $N$ is the
Ito measure of positive excursions of the linear Brownian motion: This 
is related to the fact that
the contour  process of Aldous' {\it Continuum Random Tree} is given 
by a normalized
Brownian excursion (see \cite{Al1} and \cite{Al2}), or to the 
Brownian snake construction
of superprocesses with quadratic branching mechanism (see e.g. \cite{LG1}).
For a general $\psi$, the limit theorems for the contour processes of discrete
Galton-Watson trees given 
in \cite{DuLG}, Chapter 2 
and the Ray-Knight property for the local times of $H$ proved in Theorem 1.4.1
\cite{DuLG}
both strongly justify that the height process is the right object 
to code the genealogy of (sub)critical CSBP's and that $H$ under the excursion
measure $N$ is the contour process of a continuum random tree that is called the
Lévy tree: we refer to \cite{DuLG2} for a precise definition of Lévy trees in term of
random metric $\R$-trees space (see also \cite{EvPitWin} for related topics). 
All the results about height processes used in the paper are recalled in Section 3.1.

%which relate the contour processes of Galton-Watson trees to 
%the height process $H$, and the Ray-Knight Theorem 1.4.1 \cite{DuLG} for the local times of $H$ strongly justify that the height process is the right object 
%to code the genealogy of (sub)critical CSBPs and we can build the corresponding family tree called the $\psi$-{\it Lévy tree} as a 
%random metric space which is a $\R$-tree (see \cite{DT96} 
%for a definition) and whose contour process is $H$: we refer to 
%\cite{DuLG2} for a precise discussion on this topic. 

\vspace{5mm}

\noindent
$\bullet $ The second process is a bivariate subordinator 
$(U,V)=((U_t ,V_t); \: t\geq 0)$, namely a 
$[0, \infty) \times [0, \infty)$-valued L\'evy
process started at $0$.
Its distribution is characterized 
by its Laplace exponent $\Phi$: 
$$ \E \left[ \exp \left( -pU_t -qV_t \right) \right] = \exp (-t\Phi (p,q)) .$$
$ \Phi$ is of the form 
$$ \Phi (p,q)= d p+ d'q+ \int_{(0,\infty)^2} R(dxdy)(1 -e^{-px-qy}) , $$
where $d,d' \geq 0$ and $R$ is a $\sigma$-finite measure on $(0, \infty)^2$ 
such that $\int R(dxdy) 1 \wedge (x+y) < \infty $. 
The Lévy measure $R$ plays
the role of the discrete dispatching measure $r$. 
Roughly speaking, 
think of the population as being indexed by the positive real
numbers and let us make an informal analogy with the discrete model:
If the height $t\in [0, \infty)$ in the family tree of the CSBPI
corresponds to generation $n$ in the discrete GWI-tree, then $U_t$
(resp. $V_t$) corresponds to the sum of the numbers of immigrants at
the left (resp. the right) hand of the infinite line of descent from
generation $0$ to generation $n$. Then,
a jump of $(U,V)$ occuring at time $t$ corresponds 
to a total amount $U_t-U_{t-}+V_t-V_{t-}$ of 
immigrants arriving at height $t$ in the family tree: 
$U_t-U_{t-}$ of them are put at the left hand of the infinite line of 
descent and they are the initial population of a CSBP($\psi$) ; 
$V_t-V_{t-}$ of them are put at the right hand 
of the infinite line of descent and they 
are also the initial population of an independent CSBP($\psi$). 
It implies that the real-valued subordinator $U+V$ has Laplace
exponent $\varphi$ and thus, 
$$ \varphi (p)=\Phi (p,p) .$$

  More precisely, we define $\overleftarrow{H}$
and $\overrightarrow{H}$  as follows. Let us first introduce the 
right-continuous inverses of $U$ and $V$:
\begin{equation*}
U^{-1}_t= \inf \{ s\geq 0 : \; U_s >t \} \quad 
{\rm and} \quad
 V^{-1}_t= \inf \{ s\geq 0 : \; V_s >t \} , 
\end{equation*}
with the convention $\inf \emptyset= \infty $. Let $H$ be the height process 
associated with the L\'evy process $X$ with Laplace exponent $\psi$ 
and let $(H', X')$ be an independent copy of $(H,X)$. 
We set for any $t\geq 0$, $I_t=\inf_{[0,t]}X$ and  
$I'_t=\inf_{[0,t]}X'$. Then, we define $\overleftarrow{H}$
and $\overrightarrow{H}$ by 
\begin{equation}
\label{fonccontsig}
\overleftarrow{H}_t = H_t + U^{-1 }_{-I_t } \quad 
{\rm and} \quad 
\overrightarrow{H}_t = H'_t + V^{-1 }_{-I'_t } \; , \; t\geq 0.
\end{equation}
The processes $\overleftarrow{H}$
and $\overrightarrow{H}$ are called respectively the left and the right 
height processes. The left and right height processes are continuous iff 
$U^{-1 }$ and $V^{-1 }$ are continuous, which happens iff $U$ and $V$ are 
not Poisson processes. This is equivalent to the analytical
condition
\begin{equation}
\label{UVcont}
dd'\neq 0 \quad {\rm or } \quad R\left( (0,\infty)^2\right)  =\infty .  
\end{equation}

 As for L\'evy trees, it is possible to build a $(\psi, \Phi)$-{\it immigration Lévy tree} 
via the left and right height processes:
%$\overleftarrow{H}$ and $\overrightarrow{H}$:
To each $s\in \R$ corresponds a vertex in the continuum tree at height  
$$ J_s= \un_{(-\infty , 0)} (s) \overleftarrow{H}_{-s} +  
\un_{[0, \infty )} (s) \overrightarrow{H}_s .$$
Suppose that $s\leq s'$. The common ancestor of the vertices corresponding to $s$ and 
$s'$ is situated at height $m(s, s')= \inf\{J_u ; 
\, u\in I(s,s') \} $, where $I(s, s')$ 
is taken as $[s,s']$ if $0\notin [s, s']$ and as 
$\R \setminus [s, s'] $ otherwise. Then, the distance separating 
the vertices corresponding to 
$s$ and $s'$ is given by  
$$ \dist (s, s')= J_s + J_{s'}- 2 m(s, s') .$$
Check that $\dist $ 
is a pseudo-metric on $\R$. We say that two real numbers $s$ and $s'$ 
are equivalent if they correpond to the same vertex, 
that is: $\dist (s,s')=0$. This equivalence relation is denoted by $s \thicksim s'$ 
and we formally define the 
$(\psi,\Phi)$-{\it immigration Lévy tree} as the quotient set
$ \t^*= \R / \thicksim $ equipped with the metric $\dist $ that makes it
a random Polish space. Arguing as in
\cite{DuLG2}, we can show that a.s. $ (\t^* ,\dist )$ is a real tree.

This genealogical model is clearly related to the model discussed in
\cite{Lamb} by A. Lambert where all the population is on the left hand side of the
infinite line of descent so that only one contour process is needed to
encode the family tree of the CSBPI's. In Lambert's paper, this
contour process is defined as a functional of a 
Markov process $X^*$ generalizing
the functional giving the height process introduced by Le Gall and Le Jan. 
This Markov process can be constructed either pathwise in terms
of a spectrally positive Lévy process and an independent subordinator, or
in distribution thanks to Itô's synthesis theorem.
A. Lambert also defines 
in a weak sense the local time processes of the resulting contour process
and states a generalized Ray-Knight theorem by proving they are
distributed as a CSBPI. 

\subsection{Statements of the main results.}
   In this paper, the model that we consider allows to have population 
on both sides of the infinite line of descent which turns out to be a
natural case to discuss for we can define 
continuous analogues of the discrete size-biased trees. In
particular we show
in Theorem \ref{conditionnonextinc} stated 
below that, as in the discrete case, the continuous analogues of
size-biased 
trees are the limit
of the family trees of (sub)critical CSBP's conditioned on
non-extinction. This special case strongly motivated the present
work. In the more general case, we also provide strong approximation 
for the local time processes of the right and the left contour processes
(Proposition \ref{approxtpsloc}) and we state the Ray-Knight property for them (Theorem
\ref{doubleRayKnight}). We prove functional weak convergence 
for the rescaled contour 
processes of discrete GWI-trees (Theorem \ref{convimmi}). 
We also state in Theorem \ref{strongGWI}
a general limit theorem for GWI processes (with possibly 
supercritical offspring distribution), in the same vein as 
Theorem 3.4 \cite{Gr} of Grimvall.

  More specifically, let us now 
give a detailed presentation of the main results of the
paper. The first one defines the local time processes of the two contour
processes by a strong approximation procedure. 
 \begin{proposition}
\label{approxtpsloc}
Assume that (\ref{Hcont}) and (\ref{UVcont}) hold. Then, 
there exist two jointly measurable processes 
$(\overleftarrow{L}^a_s ;\; a, s \geq 0)$ and $(\overrightarrow{L}^a_s ;
\; a, s \geq 0)$ such that: 
\begin{description}
\item
$(i)$ A.s. for any $a >0 $, $(\overleftarrow{L}^a_s ;\; 
s \geq 0)$ and 
$(\overrightarrow{L}^a_s ;\; s \geq 0)$ are continuous nondecreasing 
processes.
\item
$(ii) $ For all $T>0 $ the following limit holds in probability: 
$$ \sup_{t\leq T} \left| \epsilon^{-1}\int_0^t ds \un_{\{ a< 
\overleftarrow{H}_s \leq a+\epsilon \} }-\overleftarrow{L}^a_t 
\right|  \xrightarrow[\epsilon \rightarrow 0]{\quad } 0 .$$
The same limit holds for $\overrightarrow{H}$ and 
$\overrightarrow{L}^a$.
\item
$(iii)$ A.s. for any continuous function $g$ on $[0, \infty)$ with compact support and
for any $t\geq 0$,
\begin{equation*}
 \int_0^t ds \, g( \overleftarrow{H}_s ) = \int_0^{\infty } da \overleftarrow{L}^a_t
  g(a)  
 \quad {\rm and } \quad 
 \int_0^t ds \, g( \overrightarrow{H}_s ) = \int_0^{\infty } da
\overrightarrow{L}^a_t g(a) . 
\end{equation*}
\end{description}
\end{proposition}
We next show that the local time processes of $\overleftarrow{H}$ and 
$\overrightarrow{H}$ enjoy a ``Ray-Knight'' property: As $U$ and $V$ both
drift to infinity, so do the left and the right 
height processes and it makes sense to 
define
$$ \overleftarrow{L}^a_{\infty}=\lim_{s\rightarrow \infty} 
\overleftarrow{L}^a_{s} \quad {\rm and} \quad 
\overrightarrow{L}^a_{\infty}=\lim_{s\rightarrow \infty} 
\overrightarrow{L}^a_{s} \; ,\quad a\geq 0.$$

\noindent 
Then, the Ray-Knight theorem can be stated as follows. 

\begin{theorem}
\label{doubleRayKnight}
Assume that (\ref{Hcont}) and (\ref{UVcont}) hold. 
Then, the process 
$(\overleftarrow{L}^a_{\infty}
+\overrightarrow{L}^a_{\infty} \; ;\; a\geq 0 )$ 
is a CSBPI($\psi$,$\varphi$) started at $0$ with $\varphi (\lambda )=\Phi (\lambda ,\lambda )$ , $\lambda\geq 0$. 
\end{theorem}

%  As for L\'evy trees, it is possible to build a $(\psi, \Phi)$-{\it immigration Lévy tree} 
%via the left and right height processes:
%%$\overleftarrow{H}$ and $\overrightarrow{H}$:
%To each $s\in \R$ corresponds a vertex in the continuum tree at height  
%$$ J_s= \un_{(-\infty , 0)} (s) \overleftarrow{H}_{-s} +  
%\un_{[0, \infty )} (s) \overrightarrow{H}_s .$$
%Let $s\leq s'$. The common ancestor of the vertices corresponding to $s$ and 
%$s'$ is situated at height $m(s, s')= \inf\{J_u ; 
%\, u\in I(s,s') \} $, where $I(s, s')$ 
%is taken as $[s,s']$ if $0\notin [s, s']$ and as 
%$\R \setminus [s, s'] $ otherwise. Then, the distance separating 
%the vertices corresponding to 
%$s$ and $s'$ is given by  
%$$ \dist (s, s')= J_s + J_{s'}- 2 m(s, s') .$$
%Check that $\dist $ 
%is a pseudo-metric on $\R$. We say that two real numbers $s$ and $s'$ 
%are equivalent if they correpond to the same vertex, 
%that is: $\dist (s,s')=0$. This equivalence relation is denoted by $s \thicksim s'$ 
%and we formally define the 
%$(\psi,\Phi)$-{\it immigration Lévy tree} as the quotient set
%$ \t^*= \R / \thicksim $ equipped with the metric $\dist $ that makes him 
%a random Polish space. Arguing as in
%\cite{DuLG2}, we can show that $ (\t^* ,\dist )$ is a $\RR$-tree. 

Proposition \ref{approxtpsloc} and Theorem \ref{doubleRayKnight} are proved 
in Section 3.2 while Section 3.3 is devoted to the study of 
the continuous analogue of size
biased GW-trees,  
%that we call {\it $\psi$-size-biased 
%L\'evy trees} 
that corresponds to $(\psi, \Phi)$-immigration Lévy trees 
where $\Phi$ is given by 
$$\Phi (p,q)=\frac{\psi^*(p)-\psi^*(q)}{p-q} . $$
Here, we have set $\psi^*(\lambda)=\psi(\lambda)-\alpha\, \lambda$ and 
when $p=q$, the ratio $(\psi^*(p)-\psi^*(q))/(p-q)$ should be interpreted as
$\psi'(p)-\alpha$. So $U+V$ is
a subordinator with Laplace exponent $\varphi =\psi'-\alpha$ that is 
the immigration mechanism of the underlying CSBPI. Now, consider the height
process $H$ under its excursion measure $N$ and denote by $\zeta $ the duration of the excursion. 
As a consequence of (\ref{Hcont}) we get 
$N(\sup_{s\in[0, \zeta]} H_s >a) \in (0,\infty)$ for any $a>0$. Thus,
we can define the probability  measure  
$N_{(a)}=N(\; \cdot\;  |\sup H >a)$ (see Section 3.1 for details). 
The main result
proved in Section 3.3 can be stated as follows. 
\begin{theorem}
\label{conditionnonextinc}
Assume that (\ref{Hcont}) holds. Then, 
$$ \left( H_{t \wedge \zeta } \; ,\;  H_{(\zeta -t)_+ } \, ;\, t\geq 0 \right) 
\quad {\rm under} \quad 
N_{(a)} \; \; \xrightarrow[a \rightarrow \infty ]{\;} 
\; \; \left( \overleftarrow{H}_t \; , \; 
\overrightarrow{H}_t \, ;\, t\geq 0 \right) \; $$
weakly in $C ([0, \infty) , \R^2 ) $.
\end{theorem}
(Here $(x)_+$ stands for the non-negative part $\max (0,x)$ of $x$.) 
The proof of this theorem relies on a lemma 
(Lemma \ref{Bisbismut}) that is stated in 
Section 3.3 and that is an easy consequence of Lemma 3.4 in 
\cite{DuLG2}. Let us mention that Lemma \ref{Bisbismut} is a generalization of Bismut's decomposition of the
Brownian excursion.

   Section 4.1 and 4.2 are devoted to limit theorems 
for rescaled GWI processes and contours of GWI trees. The main result proved
in Section 4.1 is a strong invariance principle for GWI processes.  
{\bf In Section 4.1 and only in Section 4.1} 
we do not anymore restrict our attention to (sub)critical GWI-processes.
More precisely, let $(\mu_p ;p\geq 1) $ and 
$(\nu_p ;p\geq 1) $ be {\bf any} sequences of probability measures on $\N$ and let
$x\in [0, \infty)$. We denote by $(Y^{*,p}_n; n\geq 0)$, $p\geq 1$ a 
sequence of 
GWI$(\mu_p, \nu_p)$-processes started at $Y^{*,p}_0=[px]$ and we denote by 
$(\gamma_p; p\geq 1)$ an increasing sequence of positive integers.

\begin{theorem}
\label{strongGWI}
The three following assertions are equivalent:

\noindent (i) For any $t\geq 0$ the following convergence
\begin{equation}
\label{onemargi}
 p^{-1} Y^{*,p}_{[\gamma_p t]}
\build{\longrightarrow}_{p\to\infty}^{{\rm (d)}} Z^*_t ,
\end{equation}
holds in distribution in $\RR$; Here $(Z^*_t  \; ; \; t\geq 0)$ stands for a non-constant and
stochastically 
continuous process such that 
$$ \forall t>0 \; , \quad \P (Z^*_t >0)>0 \quad {\rm and } \quad \P (Z^*_t <\infty
)=1 .$$

\noindent (ii) We can find a non-constant spectrally positive L\'evy process
$X=(X_t; t\geq 0)$ with exponent $\psi$ 
and a subordinator $W=(W_t; t\geq 0)$ with exponent $\varphi $ such that 
(\ref{conservative}) holds and such that the following convergences 
\begin{equation}
\label{hypotheses}
\mu_p \left( \frac{\cdot -1}{p}\right)^{*p\gamma_p} 
\build{\longrightarrow}_{p\to\infty}^{{\rm (d)}} \P( X_1 \in \cdot \;) 
\quad {\rm and } \quad 
\nu_p \left( \frac{\cdot }{p}\right)^{*\gamma_p} 
\build{\longrightarrow}_{p\to\infty}^{{\rm (d)}} \P( W_1 \in \cdot \;)
\end{equation}
hold in distribution in $\R$ (here $*$ denotes the convolution product of measures).

\noindent (iii) There exists a non-constant and conservative 
CSBPI($\psi, \varphi$) denoted by $Y^*=(Y^*_t; t\geq 0)$,
started at $Y^*_0=x$ and such that 
\begin{equation}
\label{funcconvGWI}
\left( p^{-1} Y^{*,p}_{[\gamma_p t]}; t\geq 0\right)
\build{\longrightarrow}_{p\to\infty}^{{\rm (d)}} Y^*
\end{equation}
weakly in the cadlag functions space $\D( [0,\infty ), \R)$ endowed
with the Skorohod topology.
\end{theorem}

   In regard of Theorem 3.4 \cite{Gr} due to Grimvall that concerns limits of
GW-processes without immigration, the latter limit theorem is very
natural. However it turns out to be new. To prove $(ii) \Longrightarrow (iii)$
we adapt an argument contained in the proof of Theorem 3.4 \cite{Gr}; 
our main contribution is the proof of $(i) \Longrightarrow (ii)$.

  In Section 4.2 we prove a limit theorem for the genealogy of a 
sequence of 
(sub)critical GWI-processes: let $(\mu_p ;p\geq 1) $ be a sequence of 
offspring distributions such that $\bar{\mu}_p=\sum_{k\geq 0} k\mu_p (k)
\leq 1$ and denote by $g^{(p)}$ the corresponding generating 
functions.  Define recursively $g^{(p)}_n$ by $g^{(p)}_n=g^{(p)}_{n-1}\circ g^{(p)}$
with $g^{(p)}_0= {\rm Id}$. Let $(r_p ;p\geq 1) $ be a sequence of dispatching
distributions and let $\tau_p$ be a GWI($\mu_p, r_p$)-tree. 
For any $n\geq 0$, we also denote by $Y^{*,p}_n$ the number of non-mutants at generation $n$ in $\tau_p$. 
Recall that $(\gamma_p; p\geq 1)$ stands for an increasing sequence of positive integers. 
We suppose that 
\begin{equation}
\label{convXUV}
\mu_p \left( \frac{\cdot -1}{p}\right)^{*p\gamma_p} 
\build{\longrightarrow}_{p\to\infty}^{{\rm (d)}} \P( X_1 \in \cdot \;) 
\quad {\rm and } \quad 
 r_p \left( \frac{\cdot }{p} \; ,\;  \frac{\cdot }{p} \right)^{*\gamma_p} 
\build{\longrightarrow}_{p\to\infty}^{{\rm (d)}} \P( U_1 \in \cdot\; ; 
\; V_1 \in \cdot \;)
\end{equation}
weakly in resp. $\R$ and $\R^2$. Here $X$ stands for 
a spectrally positive L\'evy process whose exponent $\psi$ satisfies (\ref{subcri}) and (\ref{Hcont}), 
and $(U,V)$ stands for a bivariate subordinator whose exponent $\Phi$
satisfies (\ref{UVcont}). We also make the additional assumption
\begin{equation}
\label{extinctech}
\liminf_{p\to\infty} g^{(p)}_{[\delta \gamma_p]}(0)^{p}>0
\end{equation}
%We can show that (\ref{extinctech}) 
which implies that extinction times of 
$GW (\mu_p)$-processes converge in distribution in the 
$\gamma_p^{-1}$ time-scale. (\ref{extinctech})
turns out to be a necessary condition in order to have a strong convergence of
contour processes of rescaled Galton-Watson trees: see \cite{DuLG}
Chapter 2 p. 54 for a precise discussion of this point.

  Denote by $\overleftarrow{H}$ and $\overrightarrow{H}$ 
the left and the right height processes associated with a 
$(\psi, \Phi)$-immigration Lévy tree as defined by (\ref{fonccontsig})
and denote by $\overleftarrow{L}$ and $\overrightarrow{L}$ their
corresponding local times. The
main result proved in Section 4.2 is the following. 
\begin{theorem}
\label{convimmi}
Assume that (\ref{convXUV}) and  (\ref{extinctech}) hold. For any $p\geq 0$, let 
$\tau_p$ be a GW ($\mu_p$)-tree. Then, 
\begin{eqnarray*}
\left(
\left( \gamma_p^{-1} \overleftarrow{C}_{2p\gamma_p t}(\tau_p) , 
 \gamma_p^{-1} \overrightarrow{C}_{2p\gamma_pt} (\tau_p) \right)_{t\geq 0}; 
\left( p^{-1} Y^{*,p}_{[p\gamma_pa]}  \right)_{a\geq 0} \right)
\build{\longrightarrow}_{p\to\infty}^{\hspace{2mm} {\rm (d)}
\hspace{2mm}} \hspace{30mm} \\
\left(
\left(
 \overleftarrow{H}_t , \overrightarrow{H}_t
\right)_{t\geq 0} ; 
\left( \overleftarrow{L}^a_{\infty} +\overrightarrow{L}^a_{\infty} 
\right)_{a\geq 0} \right) 
\end{eqnarray*}
in distribution in $\D([0, \infty), \R^2) \times \D([0, \infty), \R)$. 
\end{theorem}
This result relies on 
combinatorial formulas stated in Section 2.2, on Theorem \ref{strongGWI}
and also on Theorem 2.3.1 \cite{DuLG} that guarantees a similar convergence for the rescaled 
contour processes of sequences of Galton-Watson trees without immigration.

  The paper is organized as follows: In Section 2.1, we set definitions and
  notations concerning discrete trees. In Section 2.2, we
discuss various codings of sin-trees that are used in the proof of Theorem 
\ref{convimmi}. In Section 3.1 we recall important properties of the 
height process that are needed to prove 
Proposition \ref{approxtpsloc}, Theorem 
\ref{doubleRayKnight} and Theorem \ref{conditionnonextinc} in Section 3.2; 
Section 3.3. Section 4.1 and Section 4.2
are devoted to the proofs of Theorem \ref{strongGWI} and Theorem 
\ref{convimmi}.

\section{Sin-trees and sin-forests}
\subsection{Definitions and examples}

Let $\N= \{ 0, 1, 2, \ldots \} $ be the set of 
the nonnegative integers, set 
$\N^* = \N \setminus \{ 0\}$ and denote the set of finite words written with 
positive integers by $\U := \{ \varnothing \} \cup \bigcup_{n\geq 1} (\N^* )^n $ . Let $u\in \U$ be the word $u_1 \ldots u_n$, $u_i \in
\N^{*}$. We denote the length of $u$ by $|u|$ : $|u|=n$.  Let $v=v_1 \ldots v_m
\in \U $. Then the word 
$uv$ stands for the concatenation of $u$ and $v$: $uv=u_1 \ldots u_n v_1
\ldots v_m $. Observe that $\U$ 
is totally ordered by the {\it  lexicographical
order} denoted by $\leq $ . A rooted ordered tree $t$ is a subset of $\U$ satisfying the following 
conditions
\begin{description}
\item{(i)} $\varnothing \in t$ and $\varnothing $ is called the {\it root} of $t$.

\item{(ii)} If $v\in t$ and if $v=uj$ for   $j\in \N^* $, then, $u\in t $. 

\item{(iii)} For every $u\in t $, there exists $k_u (t) \geq 0 $ such that $uj\in t$ for every $1\leq j\leq k_u (t)$. 
\end{description}
We denote by $\T$ the set of ordered rooted trees. 
We define on $\U$ the {\it genealogical}
order $\preccurlyeq $ by
$$\forall u,v \in \U \; , \quad u\preccurlyeq v \Longleftrightarrow \exists w\in \U \; :\; v=uw  .$$
If $u\preccurlyeq v$, we say that $u$ is an ancestor of $v$. If $u$ is distinct from the root,
it has an unique predecessor with respect to $\preccurlyeq$ who  
is called its parent and who is denoted by $\overleftarrow{u}$. 
We define the youngest common ancestor of $u$ and $v$ by the 
$\preccurlyeq$-maximal element $w\in \U$ such that $w\preccurlyeq u $ and $ w \preccurlyeq v $ and we denote it by 
$u\wedge v$. We also define the distance between $u$ and $v$ by $\dist (u,v)=|u| +|v|-2 |u\wedge v|$ and we use notation 
$\llb u,v \rrb$ for the shortest path 
between $u$ and $v$. Let $t\in \T$ and $u\in t$. We define the tree $t$
shifted at $u$ by 
$\theta_u (t)= \{ v\in \U \; :\; uv\in t\}$ and  
we denote by $[t]_u $ the tree $t$ cut at the node $u$  :  $[t]_u := \{ u\} \cup 
\{ v\in t \; :\; v\wedge u \neq u \}$. Observe that $[t]_u \in \T$. For any
$u_1, \ldots ,u_k \in t$ we also set 
$[t]_{u_1, \ldots , u_k} : =[t]_{u_1} \cap \ldots \cap [t]_{u_k}$ and
$$ [t]_n= [t]_{\{u\in t:\; |u|=n   \}}=\{ u\in t \; :\; |u| \leq n \} \; ,
\quad n\geq 0. $$

Let us denote by $\g $ the $\sigma $-field on $\T$ generated by the sets 
$\{ t\in \T :\;  u\in t\}$ , $u\in \U$. All the random 
objects introduced in this paper are defined on 
an underlying probability space denoted by $(\Omega , \f , \P)$. A
random tree is then a $\f $-$\g $ measurable mapping 
$\tau : \Omega \longrightarrow \T $.
We say that a sequence of random trees 
$(\tau_k ;\; k\geq 0 ) $ converges 
in distribution to a random tree $\tau $ iff
$$ \forall n \geq 0, \; \forall t \in \T,
 \qquad \P \left( [\tau_k ]_n = t \right) 
\xrightarrow[k\rightarrow \infty]{\quad} \P \left( [\tau ]_n = t \right) $$
and we denote it by 
$ \tau_k \xrightarrow[k\rightarrow \infty]{{\rm distr}}
\tau $.

  Let 
$\mu $ be a probability distribution on $\N$. We call {\it Galton-Watson tree} with offspring distribution 
$\mu$ (a GW($\mu $)-tree for short) any 
$\f-\g$ measurable random variable $\tau $ 
whose distribution is characterized by the following 
conditions: 
\begin{description}
\item{(i)} $\P (k_{\varnothing }(\tau )=i) = \mu (i) \; ,\; i\geq 0 $.

\item{(ii)} For every $i\geq 1$ such that $\mu (i) \neq 0$, the shifted trees $\theta_{1} (\tau) , \ldots , \theta_{i} (\tau) $  
under \\ 
$\P (\cdot \mid k_{\varnothing } (\tau )=i)$ are independent copies of 
$\tau$ under $\P$.
\end{description}

\begin{remark}
\label{indGW}
Let $u_1, \ldots , u_k \in \U$ such that $u_i \wedge u_j \notin \{ u_1, \ldots
, u_k\}$ , $1\leq i,j\leq k$, and 
let $\tau $ be a GW($\mu $)-tree. Then, as a consequence of the definition of
GW-trees, conditional on the event 
$\{ u_1, \ldots , u_k \in \tau\} $, $\theta_{u_1} (\tau) , \ldots , \theta_{u_k} (\tau) $ are i.i.d. GW($\mu $)-trees independent
of $[\tau]_{u_1, \ldots , u_k}$. 
\end{remark}

  We often consider a forest (i.e. a sequence of trees) instead of a single tree. More precisely, we define 
the forest $f$ associated with the sequence of trees $(t_l\, ; l\geq 1)$ by
the set 
$$ f= \{ (-1,\varnothing)\} \cup \bigcup_{l\geq 1} \left\{ (l,u) ,\; u\in t_l \right\} $$
and we denote by $\F $ the set of forests. Vertex $ (-1,\varnothing)
$ is viewed as a fictive root situated at generation $-1$. Let $u'=(l,u) \in f$ with $l\geq
1$; the height of $u'$ is defined by $|u'| :=|u|$ and its ancestor is defined
by $(l, \varnothing )$. For convenience, we denote it by 
$\varnothing_l:=(l, \varnothing )$. As already
specified, all the ancestors $\varnothing_1, \varnothing_2$, ...  are the 
descendants of $  (-1,\varnothing) $ and are situated at generation $0$. 
Most of the notations concerning trees extend to forests: 
The lexicographical order $\leq $ 
is defined on $f$ by taking first the individuals of $t_1$, next those of
$t_2$ ... etc and leaving $  (-1,\varnothing) $ unordered. The genealogical order $\preccurlyeq$ 
on $f$ is defined tree by tree in an obvious way. Let $v'\in f$. The youngest common ancestor of $u'$ and $v'$ 
is then defined as the $\preccurlyeq$-maximal element of $w'$ such that  $w'\preccurlyeq u'$ and $w'\preccurlyeq v'$ and 
we keep denoting it by $u'\wedge v'$. The number of children of $u'$ is $k_{u'} (f):=k_u (t_l)$ and the 
forest $f$ shifted at $u'$ is defined as the tree $\theta_{u'} (f) :=\theta_{u} (t_l)$. We also define 
$[f]_{u'}$ as the forest $\{ u'\} \cup \{ v'\in f \; :\; v'\wedge u' \neq u' \}$ and we extend in an obvious way
notations $[f]_{u_1', \dots , u_k'} $  and $[f]_n$. For
convenience of notation, we often identify
$f$ with the sequence $(t_l \, ; l\geq 1)$. When   $(t_l \, ; l\geq 1)=(t_1,
\ldots , t_k, \varnothing , \varnothing , \ldots )$, we say that $f$ is a
finite forest with $k$ elements and we abusively write $f=(t_1, \ldots ,
t_k)$.

We formally define the set of trees with a single infinite line of descent 
(called {\it sin-trees} for short) by 
$$ \T_{sin}= \left\{  t\in \T \; :\; \forall n \geq 0 , 
\quad  \# \left\{ v\in t \; :\; |v|=n \; \; {\rm and }\;  \;  \# \theta_v (t) =\infty 
\right\} =1 \; \right\} .$$
Let $t\in \T_{sin}$. For any $n\geq 0$, we denote by $u^*_n (t)$ the unique
individual $u$ on the infinite line of descent (i.e. such that $\# \theta_{u}
(t) =\infty $) situated at height $n$. Observe that
$u^*_0 (t)=\varnothing$. We use notation $\ell_{\infty} (t)=\{ u^*_n (t) ; \; n\geq
0 \}$ for the infinite 
line of descent of $t$ and we denote by $(l_n (t) \, ; n\geq 1)$ the sequence of positive integers such that
$u^*_n (t)$ is the word $l_1 (t) \ldots l_n (t) \in \U$.
We also introduce the set of {\it sin-forests} $\F_{sin}$ that is defined as the set of forests 
$f=(t_l \, ; l\geq 1)$ such that all the trees $t_l$ are finite except one
sin-tree $t_{l_0}$. We extend to sin-forests notations 
$u^*_n$, and $l_n$ by setting $l_n (f)=l_n (t_{l_0})$, $u^*_n (f)=(\, l_0, u^*_n (t_{l_0})\, )$ and $u^*_0 (f)=\varnothing_{l_0}$.

We now precisely define the Galton-Watson trees with immigration introduced in Section
1: Recall that a GWI tree is characterized by  

\noindent
$\bullet$ its {\it offspring distribution} $\mu $
on $\N$ that we suppose critical or subcritical: $\bar{\mu}=\sum_{k\geq 0} 
k\mu (k) \leq 1$;

\noindent
$\bullet$ its {\it dispatching distribution} $r$ defined on the set 
$\{ (k,j)\in \N^* \times \N^* \; :\; 1\leq j \leq k \}$ that prescribes the
distribution of the number of 
immigrants and their positions with respect to the infinite
line of descent.

\noindent
More precisely, $\tau$ is a GWI($\mu,r$)-tree if it satisfies the two following conditions:
\begin{description}
\item{(i)} The sequence $S=( \; (k_{u^*_n (\tau)} (\tau ), l_{n+1} (\tau)) \; ; \; n\geq 0)$ is i.i.d. with distribution $r$.

\item{(ii)} Conditional on $S$, the trees 
$\{ \theta_{u^*_n (\tau)i} \, (\tau), \quad   n\in \N , \;  1\leq i\leq k_{u^*_n (\tau)} (\tau ) \quad {\rm with} \quad i\neq
  l_{n+1} (\tau)\}$ are mutually independent  GW($\mu$)-trees.
\end{description}
%\noindent Informally speaking, $\tau$ can be seen as an infinite line of descent on that we have put independently 
%at each level a random number of GW($\mu$)-bushes on both sides of the infinite 
%line of descent according distribution $r$. 

We define a GWI($\mu,r$)-forest with $l\geq 1$ 
elements by the  forest $\varphi=(\tau , \tau_1 , \ldots ,\tau_{l-1}) $ 
where the $\tau_i$'s are i.i.d. GW($\mu$)-trees 
independent of the  GWI($\mu,r$)-tree $\tau$. It will be sometimes 
convenient to insert $\tau$ at random in the sequence 
$(\tau_1 , \ldots , \tau_{l-1})$ but unless otherwise specified
we choose to put the random sin-tree first in a random sin-forest.

\begin{example}
\label{sizebiased}: {\it The size-biased GW-tree.} Recall from Section 1 that 
a GW($\mu$) size-biased tree is a GWI($\mu,r$)-tree
with $r(k,j)=\mu (k)/\bar{\mu} $, $1\leq j\leq k$. 
The term ``size-biased'' can be justified by the following 
elementary result. Let $\varphi$ be a random
forest corresponding to a sequence of $l$ independent GW($\mu$)-trees and 
let  
$\varphi_{\flat}$ be a GWI($\mu,r$)-forest with $l$ elements where $r$ 
is taken as above and  
where the position of the unique random sin-tree in $\varphi_{\flat}$ is 
picked uniformly at random among the $l$
  possible choices. Check 
that for any nonnegative
measurable functional $G$ on $\F \times \U$:
\begin{equation}
\label{sizebias}
\E \left[ 
\sum_{u\in \varphi } G\left( [\varphi]_u ,  u \right)\right] = \sum_{n\geq 0} l \, \bar{\mu }^n  \E \left[ 
 G\left( [\varphi_{\flat}]_{ u^*_n (\varphi_{\flat})} , u^*_n (\varphi_{\flat})\right) \right]
\end{equation} 
and in particular 

$$ \frac {{\rm d}\P ( [\varphi_{\flat}]_n \in \cdot \, )}{{\rm d}\P (
  [\varphi]_n \in \cdot \, ) }= \frac{ Z_n (\varphi )}{  l\bar{\mu
  }^n} \; , $$ 
where $ Z_n (\varphi )=\# \{ u\in \varphi \; :\; |u|=n \} $, $ n\geq
  0$. 

\end{example}
%The words ``size-biased'' comes from the following 
%%elementary result: Let $\varphi$ be a random
%forest corresponding to a sequence of $l$ independent GW($\mu$)-trees and 
%let $\varphi_{\flat}$ be a GWI($\mu,r$)-forest with $l$ elements where $r$ is taken as above and  
%where the position of the unique random sin-tree in $\varphi_{\flat}$ 
%is picked uniformly at random among the $l$
%  possible choices. Check forest by forest  
%that for any nonnegative
%measurable functional $G$ on $\F \times \U$ 
%\begin{equation}
%\label{sizebias}
%\E \left[ 
%\sum_{u\in \varphi } G\left( [\varphi]_u ,  u \right)\right] = \sum_{n\geq 0} l \, \bar{\mu }^n  \E \left[ 
% G\left( [\varphi_{\flat}]_{ u^*_n (\varphi_{\flat})} , u^*_n (\varphi_{\flat})\right) \right]
%\end{equation} 
%and as a consequence ${\rm d}\P ( [\varphi_{\flat}]_n \in \cdot \, )/{\rm d}\P (
%[\varphi]_n \in \cdot \, ) 
%=  \# \{ u\in \varphi \; : \; |u|=n \} / l\bar{\mu }^n$. 
%\end{example}

\begin{example} 
\label{twotypes} 
: {\it A two-types GW-tree}. Let 
$\rho $ be probability measure on $\N \times \N $. 
Consider a population process with two types (say type 1 and type 2)
whose branching mechanism is described as follows : all the individuals in
the tree have the same offspring distribution; namely,   
one individual has $k$ children of type $1$ and 
$l$ children of type $2$ with probability $\rho (k, l)$. 
We order the children putting first those 
with type $1$ and next those with type $2$. 
Assume that we start with one ancestor with 
type 1. If we forget the types, 
the resulting family tree is a  GW($\mu$)-tree where $\mu $ is given by 
$$\mu (n) =\sum_{k+l=n} \rho (k,l) \; .$$
We assume that $\mu $ is (sub)critical. For any $n \geq 1$, 
denote by $A_n$ the event of a line of descent from generation 
$n$ to the ancestor that only contains individuals with type 1. 
Then we can prove easily 
$$ \tau \; \; {\rm under} \; \; \P (\, \cdot \, \mid A_n ) \quad 
\xrightarrow[n\rightarrow \infty]{{\rm distr}} \tau_{\infty} \; ,$$
where $\tau_{\infty } $ stands for a GWI($\mu , r$)-tree where $r$ is given by
$$ r(k,l) =\frac{1}{m}\sum_{j=l}^k \rho (j, k-j ) 
\quad {\rm with }\quad m = \sum_{k\geq 0} k\rho (k, \N) \; .$$
\end{example}

\begin{example} 
\label{ascendingwalk}
: {\it Ascending particle on a GW-tree}. Let 
$(\pi_n;\; n\geq 0)$ be a sequence of probability measures 
on $\N $ such that $\pi_n ( \{ 1, \ldots n \} )=1 $. Let 
$\tau $ be a critical or subcritical GW($\mu $)-tree. 
Consider a particle climbing $\tau $ at random in 
the following way: it starts at the root 
$\varnothing$; suppose it is at vertex $u\in \tau$ at time $n$, then there are 
two cases: if $k_u (\tau ) > 0$, then at time $n+1 $ the particle goes 
to $v=uj$ with probability 
$\pi_{k_u (\tau ) } (j) $; if $k_u (\tau) =0$, then 
the particle stays at $u$ 
at time $n+1 $. The particle is thus stopped at a final 
position denoted by $U$. We can show that $[\tau]_n$ conditional on
$ \{ |U|\geq n \}$ is distributed as 
$[\tau_{\infty }]_n$ where 
$\tau_{\infty} $ is a  GWI($\mu , r$)-tree with 
$$ r(k,l) = \frac{\mu (k)}{1-\mu (0)} \pi_{k} (l)\; , \; 1\leq l\leq k.$$
Consequently,
$$ \tau \; {\rm under} \; P\left( \, \cdot \, \mid  \; |U|\geq n \right) 
\quad \xrightarrow[n\rightarrow \infty]{{\rm distr}} \tau_{\infty }\; .$$
\end{example}

\subsection{Codings of sin-trees}

Let us first recall how to code a finite tree $t\in \T$. 
Let $u_0= \varnothing < u_1 < \ldots <
u_{\# t-1} $ be 
the vertices of $t$ listed in the lexicographical order. We define the {\it height process} of $t$ by $ H_n (t) = |u_n| $, 
$0\leq n<\# t$. $H(t)$ clearly characterizes the tree $t$.

  We also need to code $t$ in a third way by a path $D(t)=(D_n (t); 0\leq
n\leq \# t )$ that is defined by $  D_{n+1} (t) =
D_n (t) + k_{u_n} (t) -1$ and $D_0 (t)=0$. $D(t)$ is sometimes called
the Lukaciewicz path associated with $t$. 
It is clear that we can reconstruct $t$ 
from $D(t)$. Observe that the jumps of $D(t)$ are not smaller than
$ -1$. Moreover 
$D_n (t)\geq 0$  for any $0\leq n < \#t $ and $D_{\# t} (t)=-1$. 
We recall from \cite{LGLJ1} without proof the following formula that allows
to write the height process as a functional of $D(t)$:
\begin{equation}
\label{heightvswalk}
H_n (t)= \# \left\{ 0\leq j<n \; :\; D_j (t)= \inf_{j\leq k\leq n } D_k (t) \right\} \;  ,\; 0\leq n <\#t. 
\end{equation}

\begin{remark}
\label{loiGWheight}
If $\tau $ is a critical or subcritical GW($\mu $)-tree, then it is clear from our definition that $D(\tau)$ 
is a random walk started at $0$ that is stopped at $-1$ and whose jump distribution is given by $\rho (k)= \mu (k+1)$ , $ k\geq -1$. 
Thus (\ref{heightvswalk}) allows to write $H(\tau )$ as a functional of a 
random walk.
\end{remark}

The previous definition of $D$ and of the height process can be easily
extended to a forest $f=(t_l \, ; l\geq 1)$ of finite trees 
as follows: Since all the trees $t_l$ are finite, it is possible to list 
all the vertices of $f$ but $(-1,\varnothing)$ in the lexicographical order:
$u_0 =\varnothing_1 < u_1 < \ldots $ by visiting first $t_1$, then
$t_2$ ... etc. We then simply 
define the height process of $f$ by $H_n (f)=|u_n|$ and $D(f)$ by 
$D_{n+1} (f)=D_n(f) +k_{u_n} (f) -1$ with $D_0 (f)=0$. Set $n_p = \# t_1 +
\ldots +\# t_p $ and $n_0=0$ and observe that 
$$ H_{n_p+k} (f)= H_k
(t_{p+1}) \quad {\rm and} \quad
D_{n_p+k} (f)= D_k (t_{p+1})-p \; , \quad  
0\leq k<\# t_{p+1} \; , \, p\geq 0.$$ 
We thus see that the height 
process of $f$ is the concatenation of the height processes of the trees 
composing $f$. Moreover the $n$-th visited vertex $u_n$ is in $t_p$ iff
$p=1-\inf_{0\leq k\leq n} D_k (f)$. Then, it is 
easy to check that (\ref{heightvswalk}) remains true for every $n\geq 0$ when $H(t)$ and $D(t)$ are replaced by resp. 
$H(f)$ and $D(f)$.

  Let us now explain how to code sin-trees. 
Let $t\in \T_{sin}$. A particle visiting $t$ in the lexicographical order 
never reaches the part of $t$ at the right hand 
of the infinite line of descent. So we need two height processes or equivalently 
two contour processes to code $t$. 
More precisely, the left part of $t$ is the set $\{ u\in t:\; \exists v\in
\ell_{\infty } (t)\; {\rm s.t.} \; u\leq v\}$. This set can be listed in a
lexicographically increasing sequence of vertices denoted by
$\varnothing = u_0<u_1< \ldots $ etc.
We simply define the {\it left height process} of $t$ 
by $\overleftarrow{H}_n(t)=|u_n|$ , $n\geq 0$. $\overleftarrow{H}(t)$ completely codes the left part of $t$. 
To code the right part we consider the ``mirror
image'' $t^{\bullet}$ of $t$. More precisely, let $v\in t$ be the word
$c_1c_2 \ldots c_n $. For any 
$j\leq n$, denote by $v_j:=c_1 \ldots c_j$ the $j$-th ancestor of $v$ with 
$v_0=\varnothing$. Set 
$c^{\bullet}_j=k_{v_{j-1}} (t) -c_j +1$ and $v^{\bullet}=c^{\bullet}_1 \ldots
c^{\bullet}_n$. We then define $t^{\bullet}$ as 
$\{ v^{\bullet} , \;  v\in t \}$ and we define the {\it right height process} of $t$ as 
$ \overrightarrow{H}(t) : =\overleftarrow{H}(t^{\bullet})$.

 We next give another way to code a sin tree by two processes called  
the left contour and the right contour processes of 
the sin-tree $t$, that are  denoted
by resp. $\overleftarrow{C} (t)$ and $\overrightarrow{C} (t)$. Informally
speaking, $\overleftarrow{C} (t)$ is the distance-from-the-root process of a 
particle starting at the root and moving  
clockwise 
on $t$ viewed as a unit edge length graph embedded in the oriented half plane. We define 
$\overrightarrow{C} (t)$ as the contour 
process corresponding to the anti-clockwise journey. So we can also write 
$\overrightarrow{C} (t)=\overleftarrow{C} (t^{\bullet})$. More precisely, 
$\overleftarrow{C} (t)$ (resp. $\overrightarrow{C} (t)$) can be recovered from
$\overleftarrow{H} (t)$ (resp. $\overrightarrow{H}(t)$) through 
the following transform: Set $b_n =2n-\overleftarrow{H}_n (t)$ 
for $n\geq 0 $. Then observe that
\begin{equation}
\label{contourvsheight}
C_s (t)= \left\{   
\begin{array}{ll}
\displaystyle 
\overleftarrow{H}_n (t)-s+b_n &{\rm if} \; \displaystyle  s\in [b_n ,
b_{n+1}-1) , \; \\
\displaystyle s-b_{n+1} + \overleftarrow{H}_{n+1} (t) &{\rm if} \displaystyle
\;  s\in [b_{n+1}-1 ,b_{n+1} ] . \; 
\end{array} 
\right. 
\end{equation}

The contour process is close to the height process in the following sense: 
Define a mapping $q:\R_+\longrightarrow \Z_+$ by setting
$q(s)=n$ iff $s\in[b_n,b_{n+1})$. Check for every integer $m\geq 1$ that 
\begin{equation}
\label{contour1}
\sup_{s\in[0,m]}|\overleftarrow{C}_s(t)
-\overleftarrow{H}_{q(s)} (t)|\leq\sup_{s\in[0,b_m]}|\overleftarrow{C}_s(t)
-\overleftarrow{H}_{q(s)} (t)|\leq 1+\sup_{n\leq m}|
\overleftarrow{H}_{n+1} (t) -\overleftarrow{H}_{n} (t)|.
\end{equation}
Similarly, it follows from the definition of $b_n$ that
\begin{equation}
\label{contour2}
\sup_{s\in[0,m]}|q(s)-{s\over 2}|\leq\sup_{s\in[0,b_m]}|q(s)-{s\over 2}|
\leq {1\over 2}\sup_{n\leq m}\overleftarrow{H}_{n} (t) +1.
\end{equation}

We now give a decomposition of $\overleftarrow{H}(t)$ and 
$\overrightarrow{H}(t)$ along
$\ell_{\infty}(t)$ that is well suited to GWI-trees and that is
used in Section 3.2: Recall that $(u_n; n\geq 0)$ stands for the sequence 
of vertices of the left part of $t$ listed in the lexicographical 
order. Let us consider the set $\{ u_{n-1}^* (t)i ; 1 \leq i< l_n (t) ; n\geq 1\}$
of individuals at the left hand of $\ell_{\infty}(t)$ having a brother
on $\ell_{\infty}(t)$. To avoid trivialities, 
we assume that this set is not empty and we denote by $v_1< v_2 <
\ldots $ etc.
the (possibly finite) sequence of its elements listed in the 
lexicographical order. 

  The forest $f(t)=(\theta_{v_1} (t) , \theta_{v_2} (t) , \ldots)$ is then
composed of the bushes rooted at the left hand of $\ell_{\infty} (t)$ listed in
the lexicographical order of their roots.
Define $L_n (t):=(l_1 (t)-1) + \ldots + (l_n(t)-1)$ for any $n\geq 1$
and $L_0
(t)=0$; then, consider the $p$-th individual of 
$f(t)$ with respect to the lexicographical order on $f(t)$; 
it is easy to check that this individual is in a bush rooted in $t$ at 
height 
$$ \aal (p)=\inf \{ k\geq 0 :\; L_k(t) \geq 1- \inf_{j\leq p } D_j (f(t) )\} \; .$$ 
Thus the corresponding 
individual in $t$ is $u_{{\bf n} (p)}$ where ${\bf n} (p)$ is given by 
\begin{equation}
\label{spinaldec1}
{\bf n}(p)=p+\aal (p)
\end{equation}
(note that the first individual of $f(t)$ is labelled by $0$).  Conversely, 
let us consider $u_n$ that is the $n$-th individual of the left part of
$t$ with respect to the lexicographical order on $t$. Set 
${\bf p}(n) =\# \{ k<n : \; u_k\notin \ell_{\infty}(t)\}$ that is the number
of individuals coming before $u_n$ 
and not belonging to $\ell_{\infty}(t)$. Then
\begin{equation}
\label{spinaldec2}
{\bf p}(n)= \inf \{ p\geq 0\; :\; {\bf n}(p)\geq n \}
\end{equation}
and the desired decomposition follows:
\begin{equation}
\label{spinaldec3}
\overleftarrow{H}_n (t)=n- {\bf p}(n)  + H_{{\bf p}(n)} (\, f(t) \, ).
\end{equation}
Since $n - {\bf p} (n) = \# \{ 0\leq k <n \; : \; u_k\in \ell_{\infty } (t)
\}$, we also get 
\begin{equation}
\label{spinaldec4}
\aal ( {\bf p} (n) -1) \leq n- {\bf p} (n) \leq \aal ({\bf p} (n)).
\end{equation}
Observe that if $u_n \notin \ell_{\infty} (t)$, then 
we actually have 
$n- {\bf p} (n)=\aal ({\bf p} (n))$. The proofs of these identities follow 
from simple counting arguments and they are left to the reader . 
Similar formulas
hold for $\overrightarrow{H}(t)$ taking 
$t^{\bullet}$ instead of $t$ in (\ref{spinaldec1}), (\ref{spinaldec2}),
(\ref{spinaldec3}) and (\ref{spinaldec4}).

%\begin{figure}[ht]
%\psfrag{racine}{$0$}
%\psfrag{premier}{$1\; ,\; 0*$}
%\psfrag{deux}{$2\; ,\; 1*$}
%\psfrag{trois}{$3\; ,\; 2*$}
%\psfrag{quatre}{$4\; ,\; 3*$}
%\psfrag{cinq}{$5\; ,\; 4*$}
%\psfrag{six}{$6$}
%\psfrag{sept}{$7$}
%\psfrag{huit}{$8\; ,\; 5*$}
%\psfrag{neuf}{$9\; ,\; 6*$}
%\psfrag{dix}{$10\; ,\; 7*$}
%\psfrag{onze}{$11\; ,\; 8*$}
%\psfrag{douze}{$12$}
%\psfrag{treize}{$13\; ,\; 9*$}
%\psfrag{quatorze}{$14\; ,\; 10*$}
%\psfrag{quinze}{$15$}
%\epsfxsize=12cm
%\centerline{\epsfbox{sin.eps}}
%\caption{{\small {\it The left part of a sin-tree} $t$. The individuals 
%that are not on $\ell_{\infty} (t)$ have two labels: 
%the first one is their row in the lexicographical
%order on $t$ and the second one (taged with a star) 
%corresponds to their row in $f(t)$; individuals
%of $\ell_{\infty} (t)$ have only one label corresponding to their row in $t$.}}
%\label{leftsin.fig}
%\end{figure}

\begin{remark}
\label{spinalGWI} 
The latter decomposition is particularly useful when we consider a GWI($\mu ,r$)-tree $\tau$: In this case 
$(f(\tau), f(\tau^{\bullet}))$ is independent of $(L(\tau),
L(\tau^{\bullet}))$, $f(\tau)$ and $f(\tau^{\bullet})$
are mutually independent and 
$f(\tau)$ (resp.  $f(\tau^{\bullet})$) is a forest of i.i.d. GW($\mu
$)-trees if for a $k\geq 2$ we have $r(k,2)+\ldots 
+r(k, k)\neq 0$ (resp. $r(k,k-1)+\ldots 
+r(k, 1)\neq 0$), it is otherwise an empty forest. Moreover, the 
process $(L(\tau), L(\tau^{\bullet}))$ is a $\N \times \N$-valued random walk whose jump distribution is given by
$$ \P \left( L_{n+1} (\tau) -L_n (\tau)= m \, ; \,  L_{n+1} (\tau^{\bullet}) -L_n (\tau^{\bullet})=m' \right)= r(m+m'+1, m+1).$$
\end{remark}

\section{Continuum random sin-trees}

\subsection{The continuous time height process}

  In this section we recall from  \cite{LGLJ1} the definition of
the analogue in continuous time of the discrete 
height process defined in Section 2.2. We also recall from 
\cite{DuLG} several related 
results used in the next sections.

To define the continuous-time height process, we use an analogue of
(\ref{heightvswalk}) where the 
role 
of the random walk is played by a spectrally positive Lévy process 
$X=(X_t ; \; t\geq 0)$. The (sub)criticality of $\mu $ corresponds to the fact that $X$ does not drift to $+\infty$.
We also assume that $X$ has
a path of infinite variation (in the finite variation case, the height process is
basically a discrete process and so is the underlying 
tree: see \cite{LGLJ1} and \cite{Li} for a discussion with applications to
queuing processes). As already mentionned in the introduction, this happens if
the exponent $\psi $ of $X$ 
satisfies conditions (\ref{subcri}) 
and (\ref{Hcont}).
 By analogy with (\ref{heightvswalk}), the height process 
$H=(H_t ; \; t\geq 0)$ associated with $X$ is defined in such a way
that for
every $t\geq 0$ $H_t$ measures the size of the set: 
\begin{equation}
\label{set}
\{ s\in [0,t] \; :\;  X_{s-} =\inf_{s\leq r \leq t} X_r \} \; .
\end{equation}
To make this precise, we use a time-reversal argument: For any $t>0$, 
we define the Lévy process reversed at time $t$ by 
$$ \widehat{X}^t_s = X_t -X_{(t-s)-} \; ,\quad 0\leq s \leq t  $$
(with the convention $X_{0-}=0$). Then, $\widehat{X}^t$ is distributed as $X$ up to time $t$. Let us set for any $s\geq 0$,
$$ S_s = \sup_{r\leq s} X_r \quad {\rm and } \quad \widehat{S}^t_s = \sup_{r\leq s} \widehat{X}^t_r .$$
Then , the  set (\ref{set}) is the image of  
$$ \{ s\in [0,t] \; :\;  \widehat{S}^t_s =\widehat{X}^t_s  \} \; $$
under the time reversal operation $s\rightarrow t-s$. Recall that under our
assumptions $S-X$ is a strong Markov process for which $0$ is a regular value. So, we can
consider its local time process at $0$ that is denoted by 
$L(X)$. We define the {\it height process} by 
\begin{equation}
\label{defheight}
H_t = L_t(\widehat{X}^t ) \; ,\quad t\geq 0 .
\end{equation}
To complete the definition, we still need to specify the 
normalization of the local time: let us introduce the right-continuous inverse of $L(X)$: 
$$ L^{-1}_t =\inf \{ s\geq 0 \; :\; L_s(X) >t \}  $$
(with the convention that $\inf \emptyset = \infty $). Define $K_t$ by $X_{ L^{-1}_t} $ if $t<L_{\infty} (X)$ and by 
$\infty $ otherwise. A classical result of fluctuation theory (see \cite{Be} and \cite{Bi1}) asserts that 
$(K_t ; \; t\geq 0)$ is a subordinator whose Laplace exponent is given by 
$$ \E [ \exp ( -\lambda K_t ) ]= \exp (-ct\psi (\lambda ) /\lambda ) \; ,\quad t , \lambda \geq 0 .$$
Here, $c$ is a positive  constant that only depends on the normalization of
$L(X) $. 
We fix the normalization so that $c=1$. 
When $\beta >0$, standard results on subordinators imply for any $t\geq 0$,
$$H_t= \frac{1}{\beta} {\rm m} \left( \{ \widehat{S}^t_s \; ;\; 0\leq s \leq t\} \right) , $$
where $ {\rm m}$ stands for the Lebesgue measure on the real line.
In particular when $X$ is a Brownian motion, we see that $H$ is distributed 
as a reflected Brownian motion.

  Let us briefly recall the ``Ray-Knight theorem'' for $H$ (Theorem 1.4.2
  \cite{LGLJ1} and Theorem 1.4.1 \cite{DuLG}), that can be viewed as a
generalization of famous results about linear Brownian motion.
For any $a, t \geq 0$, we introduce the local time $L_t^a $ of $H$ at time $t$
and at level $a$ that can be defined via the following approximation
\begin{equation}
\label{approloc}
\lim_{\epsilon \rightarrow 0} \E \left[  \sup_{ 0\leq s \leq t} \left| \frac{1}{\epsilon} \int_0^s dr 
\un_{\{ a< H_r \leq a+\epsilon \}} -L_s^a \right| \right] =0 
\end{equation}
(see Proposition 1.3.3 \cite{DuLG} for details). Next, set 
for any $r\geq 0$ : $T_r=\inf \{ s\geq 0 \; :\; X_s=-r\}$ and 
$Y_a=L_{T_r}^a$ , $a\geq 0$. Then, Theorem 1.4.1 \cite{DuLG} asserts that 
$(Y_a \; ;\; a\geq 0)$  is a CSBP($\psi $)
started at $r$. 

Although the height process is in general not 
Markovian, we can still develop an excursion theory  
of $H$ away from $0$: Recall notation $I_t=\inf_{s\leq t} X_s $. Observe 
that for any $t\geq 0$, 
$H_t$ only depends on the values taken by $X-I$ on the excursion interval
that straddles $t$. Under our assumptions, $X-I$ is a strong Markov process
for which $0$ is a regular value so that $-I$ can be chosen as the local time of $X-I$
at level $0$. We denote by $N$ the corresponding excursion measure.
Let $(g_i, d_i)$, $i \in \ii$ be the excursion intervals of $X-I$ above
$0$. We can check that $\P$-a.s.
$$ \bigcup_{i\in \ii} (g_i , d_i) = \{ s\geq 0 \; :\;  X_s -I_s >0 \} =\{  s\geq 0 \; :\;  H_s >0 \} .$$
Denote by $h_i (s) =H_{g_i+s}$ , $ 0\leq s\leq \zeta_i=d_i-g_i $, $i\in \ii$
the excursions 
away from $0$. Then, each $H_i$ can be written as a functional of the
corresponding 
excursion of $X-I$ away from $0$ on $(g_i , d_i)$. 
Consequently, the point measure 
\begin{equation}
\label{Poiss}
\m (drd\omega)= \sum_{i\in \ii } \delta_{(-I_{g_i}, h_i )} (drd\omega)
\end{equation}
is a Poisson point measure with intensity $dr\, N (d\omega)$. 
Note that in the Brownian case, $ N$ is the 
Ito excursion measure of positive excursions of the reflected linear Brownian motion.

  {\it From now on until the end of this section we argue under $N$}. 
Let $\zeta $ denote the duration of the excursion. The local time
processes of the height process $(L_s^a\; ;\; 0\leq s\leq \zeta)$ , $ a\geq 0$ 
can be defined under $N$ through the same approximation as before, namely 
\begin{equation}
\label{localapprox}
\lim_{\epsilon \rightarrow 0} \, \sup_{a\geq 0} \, 
N \left(  \un_V \sup_{ 0\leq s \leq t\wedge \zeta} \left| \frac{1}{\epsilon} \int_0^s dr 
\un_{\{ a< H_r \leq a+\epsilon \}} -L_s^a \right| \right) =0 \; , \quad t\geq
0,  
\end{equation}
where $V$ is any measurable subset of excursions such that $N(V)<
\infty$. The above mentioned Ray-Knight 
Theorem for $H$ implies 
\begin{equation}
\label{localloi}
N\left( 1-\exp (-\lambda L^a_{\zeta} ) \right) =u (a, \lambda ) \; , \quad a,
\lambda \geq 0, 
\end{equation}
where we recall that $u$ is defined by (\ref{eqCSBP}). Set $ v(a) = \lim_{\lambda \rightarrow \infty} 
u (a, \lambda)$ to be a positive and finite quantity by
(\ref{Hcont}) satisfying
$a= \int_{v(a)}^{\infty} du/ \psi (u)$. By a simple argument discussed in 
Corollary 1.4.2\cite{DuLG}, we get 
\begin{equation}
\label{totheig}
N\left(  L^a_{\zeta}>0 \right) =N\left(  \sup_{s\leq \zeta} H_s >a \right) = v(a) . 
\end{equation} 
%The above mentionned ``Ray-Knight'' theorem and the limit theorem for discrete
%height processes (Theorem 2.3.1 \cite{DuLG} and Proposition 2.5.2 \cite{DuLG}) 
%suggests that $H$ is the height 
%process of a continuum random tree which can be  
%precisely defined as follows: Imagine that to any real $s\in [0,\zeta]$ 
%corresponds a vertex in the tree at height 
%$H_s$. Let $s\leq s'$. We then say that the youngest common ancestor of 
%the vertices corresponding to $s$ and $s'$ is situated at height 
%$$m(s, s')= \inf_{ s\wedge s' \leq u \leq s\vee s'} H_u $$ 
%and the distance between the vertices correponding to $s$ and $s'$ must 
%be $\dist (s, s')= H_s + H_{s'}- 2 m(s, s') $.
%We say that $s$ and $s'$ are equivalent if they correspond to the same 
%vertex in the tree: $\dist (s,s')=0$. We denote this equivalence 
%relation by $s \thicksim s'$ and define the
%continuum random tree as the quotient set 
%$ \t= \R / \thicksim $. 
%By analogy with Aldous'continuum random tree (which corresponds to the 
%quadratic branching case) we call $\t$ the $\psi$-continuum L\'evy tree
%or more simply a $\psi$-L\'evy tree.
%The pseudo-metric $\dist $ induces a metric on $\t$ that makes him 
%a (random) compact space. We can show that the metric space 
%$(\t ,\dist )$ is a $\RR$-tree: we refer to \cite{DuLG2} 
%for a precise discussion about $\RR$-trees, L\'evy trees and 
%and related topics.

   Let $a>0$ and set $N_{(a)}= N (\; \cdot \; \mid \; \sup H >a)$ that is 
a well-defined probability measure. The L\'evy tree coded by $H$ under 
$N_{(a)}$ enjoys a branching property that can be stated as follows: 
set
$$\widetilde{\tau}^a_t = \inf \{ s\geq 0 \; :\; \int_0^s dr \un_{\{ 
H_r \leq a\} }  >t \}\; .$$
We denote by $\h_a $ the $\sigma$-field generated by 
$(X_{\widetilde{\tau}^a_t } ,\; t\geq 0 )$ and by the class of the 
$N$-negligible sets of $\f$. We introduce   
the excursion intervals of $H$ above $a$:
$$ \bigcup_{ i\in \I(a)} (g_i , d_i) = \{ s\geq 0 \; :\; H_s > a \} \; $$
and we recall from Proposition 1.3.1 \cite{DuLG} the following result. 
\begin{proposition}
\label{branchprop}
(Proposition 1.3.1 \cite{DuLG}) The process 
$(L^a_s , \; s\geq 0)$ is measurable with respect to $\h_a$. 
Then, under $N_{(a)}$ and conditional on $\h_a $ the point measure 
\begin{equation}
\label{branchch}
\m_a (dld\omega ) = \sum_{i\in \I(a) } \delta_{( L^a_{g_j } , H_{(g_i + \cdot )
\wedge d_i} -a) } (dl \, d\omega)
\end{equation}
is independent of $\h_a $. Moreover it is 
a Poisson point measure with intensity
$$ \un_{[0,L^a_{\zeta}]} (l) \, dl 
\; N (d\omega) . $$ 
\end{proposition}

\begin{remark}  Proposition 1.3.1 \cite{DuLG} is actually stated
under $\P$ for the so-called {\it exploration process} 
$(\rho_t; t\geq 0)$ that is a
Markov process taking its values in 
the space of the finite measures of $[0, \infty)$ and 
that is related to the height process in the following way: 
$\P$-a.s. for any $t\geq 0$, the topological support of $\rho_t$ 
is the compact interval $[0,H_t]$. Thus it easy to deduce from 
Proposition 1.3.1 \cite{DuLG} a statement for the height process under 
$\P$ and our statement follows from the fact that 
$N_{(a)}$ is the distribution under $\P$ of the first excursion of
$H$ away from $0$ that reaches level $a$.
\end{remark}

\subsection{Proof of Theorem \ref{doubleRayKnight}.}

Recall from Section 1 the definition of the left and right height processes $\overleftarrow{H}$
and $\overrightarrow{H}$ of a ($\psi, \Phi$)-immigration Lévy tree.
We first prove Proposition \ref{approxtpsloc}.

\vspace{3mm}

\noindent
{\bf Proof of Proposition \ref{approxtpsloc}: } We only need to consider $\overleftarrow{H}$.
Recall the notation  $(g_i , d_i )$, $i\in \I$ for the excursion intervals 
of $H$ away from $0$. Set for any $a\geq 0$ and any $i\in \I$
$$\zeta_i = d_i-g_i \; ,\quad 
h_i= H_{g_i\wedge (\cdot +d_i)} \quad {\rm and} \quad  a_i= (a-U^{-1}_{-I_{g_i}})_+  , $$
where for any $x\in \R$ we have set $(x)_+=x\vee 0$. Recall notation 
$T_r = \inf \{ s\geq 0 \; : \; X_s=-r\}$ , $ r\geq 0$. Define for any $s\geq 0$,
$$\overleftarrow{L}^a_s=( -I_s-U_{a-})_+\wedge \Delta U_a \; + \; 
\sum_{i\in \I} L^{a_i}_{s\wedge d_i }-L^{a_i}_{s\wedge g_i } \; .$$
It is clear from the definition that $(s, a) \rightarrow 
\overleftarrow{L}^a_s $ is jointly measurable. Since the mappings 
$$s\longrightarrow ( -I_s-U_{a-})_+\wedge \Delta U_a \quad {\rm and} \quad  
s\rightarrow L^{a_i}_{s\wedge d_i }-L^{a_i}_{s\wedge g_i } $$ 
are non-decreasing and continuous, then $s \rightarrow \overleftarrow{L}^a_s$ is
a non-decreasing mapping and it is continuous on every open interval $(g_i,d_i)$
, $i\in \ii$ and also on $[T_{U_{a-}} ,\infty )$. Let $s\in [0,T_{U_{a-}} )$. 
Suppose that $s$ does not belong to an excursion interval of $H$ away from
$0$, that is $H_s=0$. Then, $\overleftarrow{H}_s=U^{-1}_{-I_{s}} <a$ and it
easy to check that the continuity of $H$ implies the existence of a non-empty
open interval centered around $s$ on which
$\overleftarrow{L}^a$ is a constant mapping. These observations imply $(i)$. 

%Since $[0, \infty )$ is the union of the compact 
%intervals $[g_i,d_i]$ , $i\in \ii$, we only have to check left continuity at $g_i$
%and right continuity at $d_i$ for any $i\in \ii$: to that end first observe that for any $i_0\in \ii$ 
%$$ \sup \left\{ \overleftarrow{L}^{a_i}_{d_i} \; ,\; d_i<g_{i_0}  \right\}=
%\sum_{-I_{g_i}<-I_{g_{i_0}} } L^{a_i}_{d_i }-L^{a_i}_{g_i } =
%\overleftarrow{L}^{a_{i_0}}_{g_{i_0}} .$$
%Similarly, check that 
%$$ \inf \left\{ \overleftarrow{L}^{a_i}_{g_i} \; ,\; g_i>d_{i_0}  \right\}=
%\sum_{-I_{g_i} \leq  -I_{g_{i_0}} } L^{a_i}_{d_i }-L^{a_i}_{g_i } =
%\overleftarrow{L}^{a_{i_0}}_{d_{i_0}} .$$
%Then, $(i)$ is proved. 

Point $(iii)$ follows from point $(ii)$ by standard arguments. It remains to
prove $(ii)$: By (\ref{approloc}), we see 
that for any $i\in \I$, 
$(L^{a}_{(s+g_i)\wedge d_i }-L^{a}_{g_i}\; ; \; a,s\geq 0)$
only depends on excursion $h_i$. So it makes sense to denote it by 
$(L^{a}_{s} (h_i)\; ; \; a,s\geq 0)$. 
%Then, a.s. for any $s\geq 0$ and 
%any $i\in \I$,
%$$   L^{a_i}_{s\wedge d_i }-L^{a_i}_{s\wedge g_i }= 
%L^{a_i}_{\zeta_i \wedge (s-g_i)_+} (h_i) .$$
Since $\P(U_{a-}=U_a)=1$ and by the definition of 
$\overleftarrow{H}$ we a.s. get for any $T\geq 0$, 
\begin{eqnarray*}
\sup_{0\leq t\leq T} \left| \epsilon^{-1}\int_0^t ds 
\un_{\{ a< \overleftarrow{H}_s
    \leq a+\epsilon \} }-\overleftarrow{L}^a_t \right|
 &\leq & \sum_{i\in \I}\un_{[0, U_a]} (-I_{g_i})  \\ & \times & 
\sup_{t\in [0, T\wedge \zeta_i]}  
\left| \epsilon^{-1}\int_0^t ds \un_{\{ a_i< h_i (s) \leq a_i +\epsilon \} }-
L^{a_i}_t (h_i) \right| .
\end{eqnarray*}
By conditionning on $U$, we get a.s.
\begin{equation}
\label{ineqauxi}
\EE \left[ 
\sup_{0\leq t\leq T} \left| \left. \epsilon^{-1}\int_0^t ds \un_{\{ a< \overleftarrow{H}_s
    \leq a+\epsilon \} }-\overleftarrow{L}^a_t  \right| \; \;  \right|
\;   U  \; \right] \leq 
\int_0^{U_a} dx \, {\rm n}_{\epsilon} ( a- U^{-1}_x ) ,
\end{equation}
where we have set for any $y\geq 0$,
$$ {\rm n}_{\epsilon} (y)= N \left(  \un_{V(y)} \sup_{ 0\leq s \leq T \wedge\zeta} \left| \frac{1}{\epsilon} \int_0^s dr 
\un_{\{ y< H_r \leq y+\epsilon \}} -L_s^y \right| \right), $$
with $V(y)=\{ \sup H >y\}$. By (\ref{totheig}), $N(V(y))=v(y)
<\infty$ so (\ref{localapprox}) applies and we get 
$ {\rm n}_{\epsilon} (y) \rightarrow 0$ when $\epsilon$ goes to $0$, for any fixed $y\geq 0$. 
Moreover, 
$$ {\rm n}_{\epsilon} (y) \leq  N( L^y_{\zeta }) +N\left( \epsilon^{-1} 
\int_0^{\zeta } ds \un_{\{ y< H_s \leq y+\epsilon \} } \right) =  N(L^y_{\zeta }) + \epsilon^{-1} 
\int_{y}^{y+\epsilon } da  N( L^a_{\zeta }) $$
by (\ref{localapprox}) once again. Then, use (\ref{eqCSBP}) and (\ref{localloi}) to get 
$$   N( L^a_{\zeta }) = \frac{\partial}{\partial \lambda } 
N\left( 1-e^{-\lambda L^a_{\zeta } } \right)_{|\lambda=0} =
 \frac{\partial}{\partial \lambda } 
u(a, \lambda )_{|\lambda=0} = e^{-\alpha a} \leq 1 .$$
Thus,  $ {\rm n}_{\epsilon} (y) \leq 2$ and  
%Next, deduce from the Ray-Knight theorem for $H$ recalled in Section 3.1 that
%$$ N\left( \int_0^{\zeta } ds \un_{\{ y< H_s \leq y+\epsilon \} } \right) = -\log \E \left[ 
%\exp \left( -\int_{y}^{y+\epsilon} Y_a da \right) \right] \; ,$$
%where  $(Y_a ;\; a\geq 0)$ stands for a CSBP($\psi $)  started at $Y_0 =1 $. 
%An elementary computation shows that 
%\begin{eqnarray*}
%%-\log \E \left[ \exp \left( -\int_y^{y+\epsilon } ds Y_s \right) \right] 
%&\leq & -\log \left( 1-  \int_y^{y+\epsilon } da \E [Y_a] \right) \\
%&\leq &-\log \left( 1-  \epsilon \right) ,
%\end{eqnarray*}
%since for any $a\geq 0$, 
%$$ \E [Y_a]= \frac{\partial}{\partial \lambda } 
%u(a, \lambda )_{\lambda=0} =  e^{-\alpha a} \leq 1.$$
%Then, there exists $C>0$ such that for any $\epsilon \in [0, 1/2]$, 
%$\sup_{y \geq 0} {\rm n}_{\epsilon} (y) < C$. 
$\int_0^{U_a} dx {\rm n}_{\epsilon} ( a- U^{-1}_x )$ tends a.s. 
to $0$ when $\epsilon$ goes to $0$  by dominated
convergence. $(ii)$ follows from 
(\ref{ineqauxi}) by an easy argument. 
%$(iii)$ is deduced from $(ii)$ by a
%standard argument. 
\cqfd

\vspace{3mm}

\noindent
{\bf Proof of Theorem \ref{doubleRayKnight}}: Since a.s. $(U, V)$ has no fixed discontinuity, 
we can write for any $a\geq 0 $ a.s.
\begin{equation}
\label{represent}
\overleftarrow{L}^a_{\infty} +  
\overrightarrow{L}^a_{\infty} =
\sum_{i\in \I} L^{a_i}_{\zeta_i} (h_i) + 
\sum_{j\in \I'} L^{a'_j}_{\zeta'_j} (h'_j) \; ,
\end{equation}
with an obvious notation for $h'_j$, $a'_j$ and $\zeta'_j$ , $j\in \I'$. Fix 
$0\leq b_1 <\ldots <b_n$. Deduce from (\ref{represent}) 
\begin{eqnarray}
\label{prems}
\EE \left[ \exp \left( -\lambda_1 (
\overleftarrow{L}^{b_1}_{\infty} +  \overrightarrow{L}^{b_1}_{\infty} )
-\ldots -\lambda_n  (
\overleftarrow{L}^{b_n}_{\infty} +  \overrightarrow{L}^{b_n}_{\infty} ) \right) \right]
  \\
=\EE \left[ \exp \left( -\int_0^{U_{b_n}} dx \, \omega (U^{-1}_x ) 
-\int_0^{V_{b_n}} dx \, \omega (V^{-1}_x )  \right) \right] \; , 
\end{eqnarray}
where we have set for any $s>0 $:
%\begin{eqnarray*}
$$
\omega(s) = N\left( 1- \exp \left( -\lambda_1 L^{(b_1 -s)_+ }_{\zeta }
-\ldots -\lambda_n L^{(b_n -s)_+ }_{\zeta } \right) \right). $$
% \\
%&=& -\log \E \left[ \exp \left( -\lambda_1 Y_{ (a_1 -s)_+ } -\ldots 
%-\lambda_n Y_{(a_n -s)_+ } \right) \right]  \; 
%\end{eqnarray*}
Let $(Y_a ;\; a\geq 0 )$ denote a 
CSBP($\psi $) started at $Y_0=1 $. The Ray-Knight property of the local times
of $H$ then implies 
\begin{equation}
\label{excurmultidistr}
\omega(s) =-\log \E \left[ \exp \left( -\lambda_1 Y_{ (b_1 -s)_+ } -\ldots 
-\lambda_n Y_{(b_n -s)_+ } \right) \right]  \; .
\end{equation}
\noindent
Now use the L\'evy-Ito decomposition of $(U, V)$ to get a.s.

\begin{eqnarray*}
 \int_0^{U_{b_n}} dx \omega (U^{-1}_x ) +\int_0^{V_{b_n}} dx \omega
 (V^{-1}_x ) = & & \\
(d+d') \int_0^{b_n } \omega (s) ds & + &
\sum_{0\leq s \leq b_n } ( \Delta U_s +\Delta V_s ) \omega (s) .
\end{eqnarray*}
Recall that $\varphi (\lambda )=\Phi (\lambda ,\lambda )$ , $\lambda \geq
0$ and deduce from the previous identity: 
$$ \E \left[ \exp \left( - \int_0^{U_{b_n}} dx \, \omega (U^{-1}_x ) -
    \int_0^{V_{b_n}} dx \, \omega (V^{-1}_x ) \right) \right]=
\exp \left( -\int_0^{b_n}\varphi ( \omega (s ) ) ds    \right) .$$
Denote by $(Y^*_a ;\; a\geq 0)$ a CSBPI($\psi $,$\varphi $) started at $Y^*_0
=0 $. An elementary computation (left to the reader) shows that 
$$ \E \left[ \exp \left( -\lambda_1 Y^*_{ b_1 } -\ldots 
-\lambda_n Y^*_{b_n  } \right) \right] = \exp \left( -\int_0^{b_n} \varphi ( \omega (s) )  ds \right)  \; , $$
which completes the proof of Theorem \ref{doubleRayKnight} by
(\ref{prems}) and (\ref{excurmultidistr}). \cqfd

\begin{remark}
\label{xcsbpi}
Let us explain how 
Theorem \ref{doubleRayKnight} extends to a CSBPI($\psi$,$\varphi$) 
started at an arbitrary state $r\geq0$. Set 
$\overleftarrow{H}^r_t = H_t + U^{-1 }_{(-I_t -r)_+}$. Thus, 
$\overleftarrow{H}^0=\overleftarrow{H}$. Observe that $\overleftarrow{H}^r$
coincides with $H$ up to time $T_r=\inf \{ s\geq 0: \; X_s =-r \}$. 
Then, use the Markov property at time 
$T_r$ to show that $\overleftarrow{H}^r_{T_r +\cdot}$ is independent of 
$\overleftarrow{H}^r_{\cdot \wedge T_r}$ and distributed as 
$\overleftarrow{H}^0$. Proposition \ref{approxtpsloc} and  (\ref{approloc}) make 
possible to define a local time process for 
$\overleftarrow{H}^{r}$ denoted by $(\overleftarrow{L}^{r,a}_s ; \; a,s\geq 0)$ that satisfies properties $(i)$, $(ii)$ and 
$(iii)$ of Proposition \ref{approxtpsloc}. Moreover, the previous observations imply that 
$$(\overleftarrow{L}^{r,a}_{T_r}; \; a\geq 0) \quad {\rm and} \quad 
(\overleftarrow{L}^{r,a}_{\infty} -
\overleftarrow{L}^{r,a}_{T_r}; \; a\geq 0)$$
are two independent processes: 
the first one is distributed as a CSBP($\psi$) started 
at $r$ and the the second one is a CSBPI($\psi, \varphi $) started at $0$. 
Then deduce from (\ref{eqCSBP}) and (\ref{equavv}) that the sum of these two
processes is distributed as a CSBPI($\psi, \varphi $) started at $r$.
\end{remark}

\subsection{Proof of Theorem \ref{conditionnonextinc}.}

In this section we discuss the $\psi$-size-biased 
L\'evy tree case, namely 
$$\Phi (p,q)=\frac{\psi^*(p)-\psi^*(q)}{p-q}  $$
where we have set $\psi^*(\lambda)=\psi(\lambda)-\alpha $.
%This section is devoted to a special case of 
%continuum immigration sin-trees which are the continuous analogues of size
%biased GW-trees. More precisely we call {\it $\psi$-size-biased 
%L\'evy trees} 
%the $(\psi, \Phi)$-continuum immigration sin-tree where 
%$$\Phi (p,q)=\frac{\psi^*(p)-\psi^*(q)}{p-q} $$
%whith $\psi^*(\lambda)=\psi(\lambda)-\alpha\, \lambda$.
%When $p=q$, the ratio $(\psi^*(p)-\psi^*(q))/(p-q)$ should be interpreted as
%$\psi'(p)-\alpha$, so that we see that $U+V$ is
%a subordinator with Laplace exponent $\varphi =\psi'-\alpha$ which  is 
%the immigration mechanism of the underlying CSBPI. We denote by 
%%%$\overleftarrow{H}$ and $\overrightarrow{H} $ 
%the left and the right height processes associated 
%with the $(\psi, \Phi)$-continuum immigration sin-tree with $\Phi$
%as above. The main result of the 
%section is the following theorem.
%\begin{theorem}
%\label{conditionnonextinc}
%Assume that (\ref{Hcont}) holds. Then, 
%$$ \left( H \; ,\;  H_{(\zeta -\cdot)_+ } \right) \quad {\rm under} \quad 
%N_{(a)} \xrightarrow[a \rightarrow \infty ]{\;} \left( \overleftarrow{H} , \overrightarrow{H} \right) \; $$
%where the limit holds in distribution in $C (\R_+ , \R^2 ) $.
%\end{theorem}
Let us introduce the last time under level 
$a$ for the left and the right height processes: 
$$\overleftarrow{\sigma}_a= \sup \{ s\geq 0\; : \; \overleftarrow{H}_s\leq a\}
\quad {\rm and} \quad 
\overrightarrow{\sigma}_a=\sup \{ s\geq 0\; : \; \overrightarrow{H}_s\leq
a\}. $$
One important argument in the proof of Theorem 
\ref{conditionnonextinc} is the following lemma. 
\begin{lemma}
\label{Bisbismut} 
Assume that (\ref{Hcont}) holds. Then, for any positive 
measurable function $F$ and $G$,  
$$ N\left( \int_0^{\zeta} ds F\left( 
H_{\cdot \wedge s} \right)  G \left(H_{(\zeta -\, \cdot \, ) 
\wedge (\zeta -s)} \right) \right) =
\int_0^{\infty} da e^{-\alpha a} \E \left[ F(\overleftarrow{H}_{\cdot \wedge 
\overleftarrow{\sigma}_a})G(\overrightarrow{H}_{\cdot \wedge 
\overrightarrow{\sigma}_a}) \right] \;, $$
and for any $a >0$, 
$$ N\left( \int_0^{\zeta} dL^a_s F\left( 
H_{\cdot \wedge s} \right)  G \left(H_{(\zeta -\, \cdot \, ) 
\wedge (\zeta -s)} \right) \right) =e^{-\alpha a} \E \left[ F(\overleftarrow{H}_{\cdot \wedge 
\overleftarrow{\sigma}_a})G(\overrightarrow{H}_{\cdot \wedge 
\overrightarrow{\sigma}_a}) \right] \;. $$
\end{lemma}

\noi {\bf Proof:} The second point of the lemma is an easy consequence of the
first one and of (\ref{localapprox}). Thus, we only have to prove the first
point. To that end, we introduce $M_f$ the space of all finite measures on $[0,\infty)$. If $\mu\in M_f$, we
denote by $H(\mu)\in[0,\infty]$ the
supremum  of the (topological) support of $\mu$. We also introduce a 
``killing operator''
on measures defined as follows. For every $x\geq 0$, $k_x\mu$ is the 
element of $M_f$ such that
$k_x\mu([0,t])=\mu([0,t])\wedge (\mu([0,\infty))-x)_+$ for every 
$t\geq 0$. Let $M^*_f$ stand for
the set of all measures  $\mu\in M_f$
such that $H(\mu)<\infty$ and the topological support of $\mu$
is $[0,H(\mu)]$. If $\mu\in M^*_f$, we denote by $Q_\mu$ the law 
under $\P$ of the process $H^\mu$
defined by
$$
\begin{array}{ll}
H^\mu_t=H(k_{-I_t}\mu)+H_t\;,&\quad \hbox{if }t\leq T_{\langle 
\mu,1\rangle}\;,\\
H^\mu_t=0\;,&\quad \hbox{if }t>T_{\langle \mu,1\rangle}\;,
\end{array}
$$
where $T_{\langle \mu,1\rangle}=\inf\{t\geq 0:X_t=-\langle \mu,1 \rangle\}$.
Our assumption $\mu\in M^*_f$ guarantees that $H^\mu$ has continuous 
sample paths, and
we can therefore view $Q_\mu$ as a probability measure on the space
$C_+([0,\infty))$ of nonnegative continuous functions on $[0,\infty)$.
For every $a\geq 
0$, we let $\M_a$ be the probability measure on
$(M^*_f)^2$ that is the distribution of 
$(\un_{[0,a]}(t)\,dU_t,\un_{[0,a]}(t)\,dV_t)$. 

The main argument of the proof is the key-Lemma 3.4 \cite{DuLG} 
that asserts that 
for any nonnegative measurables 
functions $F$ and $G$ on $C_+([0,\infty))$,
\begin{multline}
\label{lemmatrois}
N\Big(\int_0^\zeta ds\,F\Big( H_{(s-\cdot)_+} \Big) 
G \Big( H_{(s+\cdot )\wedge 
\zeta} \Big) \Big)\\
\quad=\int_0^\infty da\,e^{-\alpha a}\int \M_a(d\mu d\nu)\int 
Q_\mu(dh)Q_\nu(dh')F(h) G(h').
\end{multline}

We use (\ref{lemmatrois}) to complete the proof as follows: 
First observe that (\ref{lemmatrois}) implies that the height process is
reversible under $N$, namely
$$ (H_s; 0\leq s\leq \zeta) \overset{(law)}{=} 
(H_{\zeta -s}; 0\leq s\leq \zeta) \quad {\rm under} \quad N. $$
Then fix $r \geq 0$. By reversing one-by-one 
the excursions of $H$ away from $0$ on $[0, T_r]$, we get  
\begin{equation}
\label{claim}
( r+I_{(T_r - \, \cdot)_+} \; , \; H_{(T_r - \, \cdot)_+} ) 
\overset{(law)}{=}(-I_{ \cdot \wedge T_r }\; , \; H_{ \cdot \wedge T_r }) . 
\end{equation}
Next, fix  $a>0$ and set 
$\mu=\un_{[0,a]}(t)\,dU_t$. Note that for any $x$, 
\begin{equation}
\label{identi}H(k_{x} \mu) = a-U^{-1}_x \; .
\end{equation}
Deduce from (\ref{identi}), (\ref{claim}) and the definition of the left
height process that 
$$ \overleftarrow{H}_{\cdot \wedge 
\overleftarrow{\sigma}_a}  \overset{(law)}{=} H^\mu_{
(T_{\langle \mu,1\rangle} - \, \cdot )_+}  \; .   
$$
A similar identity holds for the right height process and 
the lemma follows from (\ref{lemmatrois}). \cqfd

\begin{remark} 
This lemma can be viewed as 
the continuous counterpart of identity (\ref{sizebias}). 
\end{remark}

\begin{remark} In the Brownian case $\psi (\lambda)=\lambda^2/2$, 
the left and the right height processes are two independent three-dimensional
Bessel processes and the lemma is a well-known identity due to Bismut \cite{Bis85} 
(See \cite{Ch2} for a generalization to spectrally 
L\'evy processes and \cite{Du1} to general L\'evy processes).
\end{remark}

\vspace{5mm}

\noi {\bf Proof of Theorem \ref{conditionnonextinc}:}
Let $b>0 $. For any $\omega $ in $ \D ([0, \infty ), \RR)$ we introduce 
$\tau_b (\omega ) =\inf \{ s\geq 0 :\; \omega (s) >b \} $ (with the usual 
convention $\inf \emptyset = \infty $). To simplify notations we set 
$$ \widehat{H}= H_{(\zeta - \, \cdot)_+ } \; ,  
 \;  \tau_b = \tau_b (H) \quad {\rm and} \quad 
\widehat{\tau}_b = \tau_b ( \widehat{H} ) \; .$$
We only have to prove the following convergence for any bounded
measurable function $F$, 
\begin{equation}
\label{taulim}
\lim_{a \rightarrow \infty } N_{(a)} \left( F ( H_{\cdot \wedge \tau_b } , \widehat{H}_{\cdot \wedge \widehat{\tau}_b } )  \right) = \EE \left[  F ( \overleftarrow{H}_{\cdot \wedge \overleftarrow{\tau}_b } , 
\overrightarrow{H}_{\cdot \wedge \overrightarrow{\tau}_b } ) \right] \; ,
\quad b>0 , 
\end{equation}
(with an evident notation for $\overleftarrow{\tau}_b$ and
$\overrightarrow{\tau}_b$) since it implies for any $t>0$  
$$ \lim_{b\rightarrow \infty} \lim_{ a\rightarrow \infty} 
N_{(a)}\left( \tau_b , \widehat{\tau }_b \leq t \right)
= \lim_{b\rightarrow \infty} \P \left( \overleftarrow{\tau}_b ,
  \overrightarrow{\tau}_b \leq t \right) = 0\; . $$
%which easily completes the proof of Theorem
%\ref{conditionnonextinc}. 
Let us prove (\ref{taulim}): deduce from (\ref{localapprox}) that 
$N$-a.e. the topological support of $ dL^b_{\cdot}$ is included in 
$\subset [\tau_b , \zeta -\widehat{\tau}_b
]$. Thus , by Lemma \ref{Bisbismut} 
\begin{equation}
\label{ellequ}
N \left( L^b_{\zeta } F (  H_{\cdot \wedge \tau_b } , \widehat{H}_{\cdot \wedge \widehat{\tau}_b } ) \right) =
e^{-\alpha b }  \EE \left[  F ( \overleftarrow{H}_{\cdot \wedge \overleftarrow{\tau}_b } , 
\overrightarrow{H}_{\cdot \wedge \overrightarrow{\tau}_b } ) \right]  \; .
\end{equation}
Recall from Section 3.1 the notation 
$(g_j , d_j ) $, $j\in \I(b)$ for the excursion intervals 
of $H$ above level $b$. For any $a>b$, we set
$$  Z^a_b = \# \{ j\in \I(b) :\;  \sup_{ s\in (g_j, d_j)} H_s > a \} $$
that is the number of excursions above level $b$ reaching level $a-b$. 
Deduce from Proposition \ref{branchprop} that conditional on 
$\h_b$ under $N_{(b)}$, the random variable $Z^a_b$ is 
independent of $\h_b$ and distributed as a Poisson random
variable with parameter $L^b_{\zeta }N( \sup H > a-b )=
L^b_{\zeta } v(a-b)$. Then use (\ref{ellequ}) and the obvious inclusion 
$\{ \sup H > a\} \subset \{ \sup H > b\}$ to get 
\begin{equation} 
\label{formpoiss}
 N_{(a)}\left( Z^a_b  F (  H_{\cdot \wedge \tau_b } , \widehat{H}_{\cdot \wedge \widehat{\tau}_b } ) \right) =
\frac{v(a-b)}{v(a)} e^{-\alpha b }  \E \left[   F ( \overleftarrow{H}_{\cdot \wedge \overleftarrow{\tau}_b } , 
\overrightarrow{H}_{\cdot \wedge \overrightarrow{\tau}_b } ) \right] \; .
\end{equation}
Let $C$ be a bounding constant for $F$. Then,
\begin{equation}
\label{inegalite1}
\left| N_{(a)} \left(  Z^a_b  F (  H_{\cdot \wedge \tau_b } , \widehat{H}_{\cdot \wedge \widehat{\tau}_b } ) \right)
- N_{(a)} \left(  F (  H_{\cdot \wedge \tau_b } , \widehat{H}_{\cdot \wedge \widehat{\tau}_b } ) \right) \right|
\leq C N_{(a)} \left( Z^a_b ;\; Z^a_b \geq 2 \right) \; .
\end{equation}
Now, observe that
\begin{eqnarray*}
N_{(a)} \left( Z^a_b ;\; Z^a_b \geq 2 \right) &
= \frac{v(a-b)}{v(a)} N\left( L^b_{\zeta} (1-\exp (-L^b_{\zeta} v(a-b) ) ) \right) \\
&= \frac{v(a-b)}{v(a)} (\frac{\partial }{\partial \lambda} 
u (b, 0)- \frac{\partial }{\partial \lambda}u 
(b, v(a-b) )).
\end{eqnarray*}
Since $\lim_{a\rightarrow \infty} v(a-b) =0 $ and $v(a)= u(b, v(a-b ))$, 
we get
$$ \lim_{a\rightarrow \infty} v(a-b)/v(a)=
\left( \frac{\partial }{\partial \lambda}u 
(b, 0) \right)^{-1}  
=e^{\alpha b} $$ 
and 
$$ \lim_{a\rightarrow \infty} \frac{\partial }{\partial \lambda} 
u (b, v(a-b) )=\; \frac{\partial}{\partial \lambda} u(b, 0).$$
Thus, 
$\lim_{a\rightarrow \infty} N_{(a)} \left( Z^a_b ;\; Z^a_b \geq 2 \right) =0$
and (\ref{taulim}) follows from the latter limits combined with 
(\ref{ellequ}), (\ref{formpoiss}) and 
(\ref{inegalite1}). \cqfd

\section{Limit theorems}

\subsection{ Proof of Theorem \ref{strongGWI}.}

  Recall the notations of Section 1: Let $(\mu_p ;p\geq 1) $ and 
$(\nu_p ;p\geq 1) $ be {\it any} sequences 
of probability measures on $\N$. In particular, we do not anymore assume 
that the $\mu_p$'s  are (sub)critical. Let 
$(\gamma_p; p\geq 1)$ be an increasing sequence of positive integers.
Denote by $g^{(p)}$ and $f^{(p)}$ the generating functions of resp. $\mu_p$
and $\nu_p$. Let $x\in [0, \infty)$. Recall that for any $p\geq 1$, 
$(Y^{*,p}_n; n\geq 0)$ stands for a 
GWI$(\mu_p, \nu_p)$-process started at $Y^{*,p}_0=[px]$. We also 
need to introduce for any $p\geq 1$ a random walk $(W^{p}_n; n\geq 0)$ 
independent of the $Y^{*,p}$'s,  
started at $0$ and whose jumps distribution is $\nu_p $. We denote by 
$(Y^{p}_n; n\geq 0)$, $p\geq 1$ a sequence of 
GW$(\mu_p)$-processes started at $Y^{p}_0=p$.

   One important ingredient of the proof is Theorem 3.4 \cite{Gr} due to
Grimvall that is the exact analogue of Theorem \ref{strongGWI} without
immigration. For convenience of notation we re-state it as a lemma. 
\begin{lemma}
\label{Grimvall}(Theorem 3.4 \cite{Gr}) 
The three following assertions are equivalent 

\noindent (a) For any $t\geq 0$, 
\begin{equation}
\label{onemargin}
 p^{-1} Y^{p}_{[\gamma_p t]}
\build{\longrightarrow}_{p\to\infty}^{{\rm (d)}} Z_t
\end{equation}
where the process $(Z_t  \; ; \; t\geq 0)$ 
is a stochastically continuous process such that 
$$ \forall t>0 \; , \quad \P (Z_t >0)>0 \quad {\rm and } \quad \P (Z_t <\infty
)=1 .$$
\noindent (b) There exists a spectrally positive L\'evy process
$X=(X_t; t\geq 0)$ with exponent $\psi$ 
satisfying (\ref{conservative}) such that the following convergence 
\begin{equation}
\label{hypothe}
\mu_p \left( \frac{\cdot -1}{p}\right)^{*p\gamma_p} 
\build{\longrightarrow}_{p\to\infty}^{{\rm (d)}} \P( X_1 \in \cdot \;) 
\end{equation}
holds weakly in $\R$.

\noindent (c) There exists a conservative stochastically continuous 
CSBP($\psi $) denoted by $Y=(Y_t; t\geq 0)$
started at $Y_0=1$ such that the convergence 
\begin{equation}
\label{funcconvGW}
\left( p^{-1} Y^{p}_{[\gamma_p t]}; t\geq 0\right)
\build{\longrightarrow}_{p\to\infty}^{{\rm (d)}} Y
\end{equation}
holds weakly in  $\D( [0,\infty ), \R)$. 
\end{lemma}

\begin{remark}
Theorem 3.4 \cite{Gr} is stated with a different scaling: we refer to the proof
Theorem 2.1.1 \cite{DuLG} to derive Lemma \ref{Grimvall} from Theorem 3.4
\cite{Gr}. 
\end{remark}

\vspace{3mm}

\noindent 
Since obviously Theorem \ref{strongGWI} 
$(iii)\Longrightarrow$ Theorem \ref{strongGWI}$(i)$, we only have to prove 
$(i)\Longrightarrow (ii)$ and $(ii)\Longrightarrow (iii)$.

\vspace{3mm}

\noindent
{\bf Proof of $(i)\Longrightarrow (ii) $}: 
For any $t, \lambda \in [0, \infty)$ and any $p\geq 1$, we set 
$$ u_p(t ,\lambda) =-p\log \left( g^{(p)}_{[\gamma_p t]}(
e^{-\lambda /p} )\right) \; , \quad \varphi_p(\lambda)=-\gamma_p  
\log \left(f^{(p)} (e^{-\lambda /p})\right) $$
 $$ b_p(t ,\lambda) = \int_0^{[\gamma_p t]/ \gamma_p} ds \; 
\varphi_p (  u_p (s ,\lambda) ) \;, \quad  
d_p(t ,\lambda) =-\log \E \left[ 
e^{-\lambda p^{-1}Y^{*,p}_{[\gamma_p t]} }  \right] $$
and $d (t ,\lambda) =-\log \E [ 
\exp (-\lambda Z^*_t )  ]$. First deduce from (\ref{transGWI})
\begin{equation}
\label{1analy}
\frac{[px]}{p} u_p(t ,\lambda ) + b_p(t ,\lambda)=d_p(t ,\lambda) .
\end{equation}
The convergence of Theorem \ref{strongGWI} $(i)$ combined with Dini's theorem 
implies the following assertions:

\noindent
$\bullet$ For any $\lambda \geq 0$, $d(0,\lambda ) =x\lambda$ 
and $d(\cdot ,\lambda)$ is continuous on $[0, \infty)$.

\noindent
$\bullet$ For any $t,\lambda >0$, $d(t,\lambda)\in (0,\infty)$ and
$\lim_{\lambda \rightarrow 0} d(t,\lambda)=0$.

\noindent
$\bullet$  For any $t \geq 0$, $d_p(t,\cdot)
\build{\rightarrow}_{p\to\infty}^{\,}d(t,\cdot)$ uniformly on every compact
subsets of $[0, \infty)$.

  Let $T$ be any denumerable dense subset of $(0, \infty)$ and let $E$ be any infinite
subset of $\N$. By use of Helly's selection
theorem combined with Cantor's diagonal procedure, we can find an increasing sequence 
$(p_k; k\geq 1)$ of elements of $E$, a set of measures $({\rm m_t}, t\in T)$ on 
$[0, \infty )$ and a measure ${\rm n}$ on $[0, \infty )$ such that 
${\rm m_t}([0, \infty )) \leq 1$, $t\in T$,   
${\rm n}([0, \infty )) \leq 1$ and such that for all $t\in T$,
\begin{equation}
\label{Helly1}
\forall r\in [0, \infty ) \; \; {\rm s. t.} \; \; {\rm m_t} (\{ r \})=0
\; , \; 
\lim_{k\rightarrow \infty} \P ( p_k^{-1} Y^{p_k}_{[\gamma_{p_k} t]} \leq r)
=  {\rm m_t} ([0, r] ) \; ,
\end{equation}
\begin{equation}
\label{Helly2}
\forall r\in [0, \infty ) \; \; {\rm s. t.} \; \; {\rm n} (\{ r \})=0
\; , \; 
\lim_{k\rightarrow \infty} \P ( p_k^{-1} W^{p_k}_{\gamma_{p_k}} \leq r)
=  {\rm n} ([0, r] ). 
\end{equation}
Define for any $\lambda \geq 0$ and any $t\in T$, 
$$  u(t,\lambda)= -\log \int_{[0, \infty)} e^{-\lambda y}{\rm m_t}(dy) 
\quad {\rm and} \quad \varphi (\lambda)= -\log \int_{[0, \infty)} e^{-\lambda y}
{\rm n}(dy) ,$$
with the convention $-\log (0)=\infty $ so that $m_t=0$
iff $u(t,\lambda)=\infty$ for a certain $\lambda \geq 0$. 
By Dini's theorem and standard monotonicity arguments, we deduce from 
(\ref{Helly1}) and (\ref{Helly2}) that for any $t\in T$ the following
convergences hold 
\begin{equation}
\label{Helly3}
u_{p_k}(t , \cdot ) 
\build{\longrightarrow}_{k\to\infty}^{\; } u(t , \cdot)
\quad {\rm and} \quad  \varphi_{p_k} 
\build{\longrightarrow}_{k\to\infty}^{\; } \varphi
\end{equation}
as $[0,\infty]$-valued functions 
uniformly on every compact subsets of the open interval $(0, \infty)$. 

\vspace{3mm}

We consider two cases: $x\neq 0$ and $x=0$ and {\bf we first suppose $x\neq
  0$}. By (\ref{1analy}), we get 
$$\frac{[px]}{p} u_p(t ,\lambda ) \leq d_p(t ,\lambda) .$$
Then we pass 
to the limit along $(p_k; k\geq 1)$ to show that $u(t,\lambda)\leq x^{-1}
d(t ,\lambda)<\infty$ , $t\in T$ , 
$\lambda >0$. Thus ${\rm m_t}\neq 0$ for any $t\in T$ 
and it makes sense to define 
the function $b$ on $T\cup \{ 0\} \times [0, \infty)$ by 
$$  b(t , \lambda )= d(t , \lambda )-xu(t , \lambda ) \quad {\rm if}
\quad t\in T,
\lambda \geq 0 $$
and $ b(0, \lambda )=0$, $\lambda \geq 0$. Deduce from (\ref{Helly3})
that for any $t\in T$ 
\begin{equation}
\label{Hellyly}
b_{p_k}(t , \cdot ) 
\build{\longrightarrow}_{k\to\infty}^{\; } b(t , \cdot) \; ,
\end{equation}
uniformly on every compact subsets of the open interval $(0, \infty)$. We first prove the following claim.

\vspace{3mm}

{\bf Claim 1}: $\quad \quad$ There exists $t_0 \in T$ such that ${\rm m_{t_0}}\neq \delta_0$ 

\vspace{3mm}

\noindent
{\bf Proof of Claim 1}: Suppose that $u(t,\lambda)=0$ , for any $t\in T$ and
any $\lambda >0$.  Denote by 
$q_p\in [0,1]$ the smallest solution in $[0,1]$ of 
$g^{(p)} (z)=z$. Observe that 
$t\rightarrow  g^{(p)}_{[\gamma_p t ]}(z)$ is
non-decreasing for $0\leq z\leq q_p$ and non-increasing for 
$q_p\leq z\leq 1$. Thus, for any $p\geq 1$ and any $\lambda \geq 0$, 
$u_p(\cdot , \lambda)$ is monotone. Since $T$ is dense, then a standard monotonicity
argument implies that 
\begin{equation}
\label{uunif}
u_{p_k}(t , \lambda ) 
\build{\longrightarrow}_{k\to\infty}^{\; } 0 \quad , \; t\geq 0, \;
\lambda >0 \;. 
\end{equation}
We take $u(\cdot , \lambda)=0$, $\lambda \geq 0$, then we also get 
\begin{equation}
\label{egalili}
b(t,\lambda)=d(t,\lambda) \; ,\quad  t\geq 0 \, ,\;  \lambda >0
\end{equation}\
and since the $b_{p_k}(\cdot , \lambda )$'s are non-decreasing and 
$d(\cdot , \lambda )$ is continuous, 
(\ref{Hellyly}) holds for any $t\geq 0$.

   Now, set $s_p= \gamma_p^{-1}([\gamma_p (t+s)]-[\gamma_p t])$ and
 use the Markov property for $Y^{*,p}$ at time $[\gamma_p t]$ to get 
\begin{equation}
\label{Markovimmdis}
 d_p(s +t,\lambda)=d_p(t, u_p(s_p,\lambda))+ b_p(s_p,\lambda).
\end{equation}
Since $s_{p_k}\rightarrow s$, since the $b_{p_k}(\cdot , \lambda )$'s
and the  $u_{p_k}(\cdot , \lambda )$'s are monotone and since their
limits are continuous, we get 
$$ u_{p_k}(s_{p_k} , \lambda ) 
\build{\longrightarrow}_{k\to\infty}^{\; } 0 \quad {\rm and} \quad 
b_{p_k}(s_{p_k} , \lambda ) 
\build{\longrightarrow}_{k\to\infty}^{\; } d(s , \lambda) , 
\quad  t\geq 0  ,\;  \lambda >0. $$
Use this to pass to the limit in (\ref{Markovimmdis}) to get 
$$ d(s+t,\lambda)=d(s,0)+ b(t,\lambda)= d(t,\lambda) , \quad t,s\geq 0
, \, \lambda >0 .$$
It then implies that $d(s,\lambda)=
x\lambda $, $s \geq 0$ and the process $Z^*$ has to be a constant process 
which contradicts 
the assumptions of Theorem \ref{strongGWI} $(i)$. \cqfd

%
% and then $b(\cdot ,\lambda)$ is non-decreasing and continuous. 
%Use monotonicity
%again to get for any 
%Then for any $t\in T$
%\begin{equation}
%\label{bmono}
% b_{p_k}(t, \cdot ) 
%\build{\longrightarrow}_{k\to\infty}^{\; } d(t, \cdot)
%\end{equation}
%uniformly on every compact subsets of the open interval 
%$(0, \infty)$. 
%Let $t\in T$ and
%$s\in [0, \infty)$. For any $p\geq 1$, we may and will choose  $s_p$ such that $|s_p -s| \leq 1/ \gamma_p$ and $[\gamma_p (s_p +t)]
%=[\gamma_ps]+ [\gamma_p t]$. Use the Markov property for $Y^{*,p}$ at 
%$[\gamma_p s]$ to get
%\begin{equation}
%\label{Markovimmdis}
% d_p(s_p +t,\lambda)=d_p(s, u_p(t,\lambda))+ b_p(t,\lambda).
%\end{equation}
%Since $d(\cdot ,\lambda)$ is continuous and since the $ b_{p_k}(\cdot , \lambda )$'s
%are non-decreasing, (\ref{Hellyly}) implies that $d(\cdot ,\lambda)$ is 
%non-decreasing. By a simple monotonicity argument, 
%$ d_{p_k}(s_{p_k} +t,\lambda)$ tends to $d(s+t,\lambda)$ when 
%$k$ goes to $\infty$. By (\ref{Hellyly}) and (\ref{egalili}), we get 
%$$ d(s+t,\lambda)=d(s,0)+ b(t,\lambda)= d(t,\lambda) \;, \quad t\in T \, , \,s\geq 0 \, , \, \lambda >0 .$$
%Since $T$ is dense and $d(\cdot , \lambda)$ is continous, we 
%deduce that $d(s,\lambda)=
%x\lambda $, $s \geq 0$ and $Z^*$ has 
%to be a constant process which contradicts 
%the assumptions of Theorem \ref{strongGWI} $(i)$. \cqfd

\vspace{3mm}

{\bf Claim 2}: $\quad \quad$ ${\rm n}\neq 0.$
\vspace{3mm}

\noindent
{\bf Proof of Claim 2}: Recall from the proof of {\it Claim 1} that 
for any $p\geq 1$ and any $\lambda \geq 0$, 
$u_p(\cdot , \lambda)$ is monotone.
Since the $f^{(p)}$'s are non-decreasing, we get for any 
$p\geq 1$ and any $t\geq 0$, 
\begin{equation}
\label{cmonoo}
\frac{[t\gamma_p]}{\gamma_p} \varphi_p \left( \lambda \wedge u_p( t, \lambda) \right) \leq 
b_p( t, \lambda) \leq d_p( t, \lambda)
\end{equation}
(use (\ref{1analy})) for the right member). 
Let $t_0 \in T$ satisfying {\it Claim 1}. Then, for any $\lambda >0$, 
$u(t_0, \lambda )\in (0, \infty)$. So it makes sense to pass to the limit 
in (\ref{cmonoo}) along $(p_k;k\geq 1)$ with $t=t_0$. We obtain
\begin{equation}
\label{cmono}
t_0  \, \varphi  \left( \lambda \wedge u( t_0, \lambda) \right) 
\leq d ( t_0, \lambda) <\infty, 
\end{equation}
which implies the claim. \cqfd

\vspace{3mm}

%\noindent
%Define the function $b(\cdot , \lambda )$ on $T\cup \{ 0\} \times (0, \infty)$ by 
%$$  b(t , \lambda )= d(t , \lambda )-xu(t , \lambda ) $$
%and $ b(0, \lambda )=0$. 
%\vspace{3mm}

{\bf Claim 3}: 
For any $\lambda >0$, $b(\cdot , \lambda )$ extends to a 
non-decreasing continuous function on $[0, \infty )$. 

\vspace{3mm}

\noindent 
Provided that {\it Claim 3} holds, Dini's theorem combined with a monotonicity
argument implies the following
convergence :
\begin{equation}
\label{Helly333}
b_{p_k}(t,\cdot) 
\build{\longrightarrow}_{k\to\infty}^{\; }
b(t,\cdot) \; \quad t\in [0, \infty)
\end{equation}
holds uniformly on every compact subsets of the open interval $ (0, \infty)$.

\vspace{3mm}

\noindent
{\bf Proof of Claim 3}: 
%Then, (\ref{1analy}) implies
%that ${\rm n}, {\rm m_t}\neq 0$ and either 
%${\rm n}\neq \delta_0$ either ${\rm m_t}\neq \delta_0$ , $t\in T$. 
%Consequently 
%$a$ and $c$ are finite and either $a$ either $c$ is 
%positive. Then we can define  
%$b(t,\lambda)= -\log \E[ \exp(-\lambda Z^*_t )] -u(t,\lambda)$ for any 
%$t\in T$ and $\lambda \geq 0$. Fix $\lambda \in (0,\infty)$. 
%We claim that 
%\begin{equation}
%\label{Helly3}
%\lim_{k\rightarrow \infty} b_{p_k}(t,\lambda)=b(t,\lambda) \; ,
%\quad t\in T , \lambda \geq 0.
%\end{equation}
%{\it Proof of the claim}: 
Fix $p\geq 1$ and $\lambda >0$ and recall that $u_p(\cdot , \lambda)$ is
monotone and that $\varphi_p$ is non-decreasing. 
Thus, we get for any $0 \leq s<t$, 
\begin{eqnarray*}
 0\leq b_p(t , \lambda)- b_p(s , \lambda)&\leq &
\frac{[\gamma_pt]-[\gamma_ps]}{\gamma_p} \varphi_p ( u_p(t , \lambda)\vee 
 u_p(s , \lambda)) \\
&\leq & \frac{[\gamma_pt]-[\gamma_ps]}{\gamma_p}
 \varphi_p \left( \frac{p}{[px]}  (d_p(t , \lambda)\vee 
 d_p(s , \lambda) )\right), 
\end{eqnarray*}
by (\ref{1analy}). Next, by (\ref{Hellyly}) 
and {\it Claim 2}
%holds for any $t\in T$. 
%Thus, we can pass to the limit along $(p_k;k\geq 1)$ in the previous
%inequality to get 
\begin{equation}
\label{sandwich}
0\leq b(t , \lambda)- b(s , \lambda) \leq (t-s) \, 
\varphi \left( x^{-1} (d(t , \lambda)\vee 
 d(s , \lambda) )\right) \; ,\quad \lambda >0 , \; s,t\in T, 
\end{equation}
%by (\ref{Helly3}) and since $\varphi$ is finite by {\it Claim 2}. This inequality
which completes the proof of the claim. \cqfd

\vspace{3mm}

Thus, it makes sense to extend the definition of $u$ on $[0, \infty )\times [0, \infty )$ by setting 
\begin{equation}
\label{aaaaa}
u(t,\lambda ):= x^{-1}\left( d(t , \lambda)-b(t , \lambda) \right) \; \quad
t\geq 0 , \;\lambda >0 
\end{equation}
and $u(t, 0)=0$, $t\geq 0$. Observe that 
$b(t , \lambda)\leq d(t , \lambda)$ and $u(t , \lambda)\leq x^{-1} d(t ,
\lambda)$, for any $t\geq 0$ and any $\lambda >0$. If $\lambda $
goes to $0$, then 
$u(t,0+)=b(t,0+)=0$ for any $t\geq 0$. It implies in particular that 
for any $t\in T$ ${\rm m}_t$ is a
probability distribution. 
This result combined with (\ref{Helly333}) also implies that the
convergence  
\begin{equation}
\label{Helly4}
u_{p_k}(t,\cdot) 
\build{\longrightarrow}_{k\to\infty}^{\; }
u(t,\cdot) \; ,\quad t\in [0, \infty)
\end{equation}
holds uniformly on every compact subsets of the closed interval 
$ [0, \infty)$. Now, set for any $t,y\in
[0, \infty)$,  
$$ Y_t(k,y)=p_k^{-1}Y^{p_k}_{[\gamma_{p_k} t]}  \quad {\rm with} \quad
Y^{p_k}_0=[p_k y] . $$
As a consequence of (\ref{Helly4}) and of the continuity theorem for Laplace
exponents (see \cite{Fe} p. 431), there exists a family of probability
measures $(P_t (y,
dz); t,y\geq 0)$ on $[0, \infty)$ 
%such that $\int P_t (y,dz)\exp (-\lambda z)
%=\exp (-yu(t,\lambda ))$ and 
such that the distribution of $  Y_t(k,y)$
converges weakly to $P_t (y,dz)$. 
In particular, 
$ P_t (1,dz)={\rm m}_t (dz)$ , $t\in T$. 
%and $\int P_t (y,dz)\exp (-\lambda z)
%%=\exp (-yu(t,\lambda ))$. 
Since $(t,y)\rightarrow \int P_t (y,dz)\exp
(-\lambda z) =\exp (-yu(t,\lambda ))$ is continuous for any $\lambda \geq 0$, the mapping $(t,y)\rightarrow
P_t (y,B)$ is measurable for any Borel set $B\subset [0, \infty)$.
Moreover, the Markov property and the branching
property for the $Y^{p_k}$'s imply that for any $ t,s,y,y'\geq 0 $:
\begin{equation}
\label{Helly5}
\int P_t(y, dy') P_s (y', dz)=P_{t+s}(y, dz) \quad {\rm and} \quad P_t (y,dz) * P_t (y',dz)=P_t (y+y',dz).
\end{equation}
By Theorem 4 \cite{Sil68} of Silverstein  
(see also the correspondence between spectrally positive Lévy processes and
CSBPs in Theorems 1 and 2 \cite{La2}) there exists a
spectrally positive Lévy process $X$ with exponent $\psi$ satisfying
(\ref{conservative}) such that $u (t, \lambda )$ is the unique nonnegative 
solution of the differential equation (\ref{equavv}). Then 
the $(P_t (y,dz); t,y\geq 0)$ are the 
transition kernels of a non-zero conservative CSBP($\psi $). The
Lévy-Khintchine formula implies that $\psi$ is of the form
\begin{equation}
\label{levKh}
\psi (\lambda) = \alpha_0 \lambda+ \beta \lambda ^2 \int_{(0,\infty)} \pi (dr) \left(
e^{-\lambda r} -1+ \lambda r\un_{\{ r<1\}} \right)
\end{equation}
with $\alpha_0 \in \R$, $\beta \geq 0$ and $\int  \pi (dr) 1\wedge r^2
<\infty$. 
%\begin{equation}
%\label{equaa}
%\frac{\partial}{\partial t}  u(t,\lambda ) = -\psi 
%\left( a (t, \lambda ) \right)
%\quad {\rm and } \quad  a (0,\lambda ) =\lambda \; ,\; \lambda\geq 0.
%\end{equation}

 Concerning the immigration exponent, 
use (\ref{cmono}) to get $ \varphi (0+)=d(t,0+)=0$. Thus, ${\rm n}$ is a true
probability distribution that has to be infinitely divisible on $[0, \infty)$
since it is obtained as 
a weak limit of marginals of rescaled random walks. $\varphi$ is therefore the
Laplace exponent of a conservative 
subordinator denoted by $W$. It has to be of the form
\begin{equation}
\label{levKhfi}
\varphi (\lambda) = \kappa\lambda + \int_{(0,\infty)} \rho (dr) \left( 1-
e^{-\lambda r}\right)
\end{equation}
with $\kappa\geq 0$ and $\int  \rho (dr) 1\wedge r <\infty$. Deduce from (\ref{1analy}), (\ref{Helly3}) and
(\ref{Helly4}) that 
\begin{equation}
\label{eqdd}
d(t, \lambda)=x u(t,\lambda) + \int_0^t \, ds \, \varphi \left( u(s,\lambda) \right)
\; , \quad \lambda , t\geq 0. 
\end{equation}

   We now need to show {\bf uniqueness} for the limiting functions $u$, $\psi$ and $\varphi$: 
let $\widetilde{u}$, $\widetilde{\psi}$ and $\widetilde{\varphi}$ be obtained by
repeating the previous procedure from another denumerable dense subset $\widetilde{T} \subset [0, \infty)$ and another
subsequence $(\widetilde{p_k};k\geq 1)$. Thus, we must have
$$ d(t, \lambda)=x u(t,\lambda) + \int_0^t ds \, \varphi \left( u(s,\lambda) \right)=
x \widetilde{u}(t,\lambda) + \int_0^t ds \, 
\widetilde{\varphi}\left(\widetilde{u} (s,\lambda) \right) \; ,\quad t, \lambda \geq 0.$$
Differentiate twice the latter equation at $t=0$ to get 
$$   -x\psi+\varphi =-x\widetilde{\psi}+\widetilde{\varphi}
\quad {\rm and } \quad
( x \psi'-\varphi') \psi =(
x\widetilde{\psi}'-\widetilde{\varphi}')\widetilde{\psi} .$$
Differentiate the first expression and deduce from the second one the 
following equation:
\begin{equation}
\label{identiti}
( x \psi'-\varphi') (\psi -\widetilde{\psi} )=0.
\end{equation}
Suppose that $x \psi'=\varphi'$ on a non-empty open interval $(a,b)$. 
Differentiate twice this expression to get 
$$ \int_{(0,\infty)} \rho (dr)r^3 e^{-\lambda r} =-x\int_{(0,\infty)} \pi
(dr)r^3 e^{-\lambda r} \; ,\quad \lambda \geq 0, $$
by (\ref{levKh}) and (\ref{levKhfi}). Thus, $\rho=\pi=0$, $\beta =0$
and $x\alpha_0=\kappa$, which imply that $Z^{*}$ is
constant: This contradicts the assumption of Theorem \ref{strongGWI}
$(i)$. Accordingly,  $x \psi'$ and $\varphi'$ must differ at a point
and by continuity $x \psi'\neq \varphi'$ on a non-empty open
interval. Thus by (\ref{identiti}), we get $\psi =\widetilde{\psi}$ 
and $\widetilde{\varphi}=\varphi$, which implies the desired uniqueness.
\cqfd

\vspace{3mm}
%If we suppose that $x \psi'=\varphi'$, then $x \psi-\varphi=0$ and the process
%$(X_{xt}; t\geq 0)$ has to be distributed as the process $(-W)$ and this is
%not possible. Thus, $x
%\psi'\neq \varphi'$. A similar argument implies $x\widetilde{\psi}'\neq
%\widetilde{\varphi}'$ and therefore $\psi =\widetilde{\psi}$ and
%$\varphi = \widetilde{\varphi}$ by an easy argument. The proof of the uniqueness is then achieved. 
%\cqfd
%\vspace{3mm}

  Thus, we have shown in the $x\neq 0$ case that there exist 
a non-constant CSBP($\psi$): $Y=(Y_t;t\geq 0)$ started at $Y_0=1$ and a 
subordinator $W$ with Laplace exponent $\varphi$ both satisfying
(\ref{conservative}) such that for any denumerable dense subset 
$T\subset (0,\infty)$ and any
infinite subset $E\subset \N$, we can find an increasing sequence 
$(p_k; k\geq 1)$ of elements of $E$ (depending on $T$) that satisfies for any
$t\in T$ and any $\lambda \geq 0$: 
$$\lim_{k\rightarrow \infty}
\E \left[ e^{-\lambda p_k^{-1} Y^{p_k}_{[\gamma_{p_k} t]} } \right] = \E
\left[ e^{-\lambda Y_t }\right] \quad {\rm and} \quad
\lim_{k\rightarrow \infty}
\E \left[ e^{-\lambda p_k^{-1} W^{p_k}_{\gamma_{p_k}}} \right] = \E
\left[ e^{-\lambda W_1 }\right] .$$ 
It easy to prove that these limits imply  for any $t\geq 0$ that the
following convergences 
$$ p^{-1} Y^{p}_{[\gamma_p t]}
\build{\longrightarrow}_{p\to\infty}^{{\rm (d)}} Y_t \quad {\rm and } \quad 
 p^{-1} W^{p}_{\gamma_p }
\build{\longrightarrow}_{p\to\infty}^{{\rm (d)}} W_1 $$
hold in distribution in $\R$. Theorem \ref{strongGWI} $(ii)$ follows from Lemma \ref{Grimvall} $(a)\Longrightarrow
(b)$. \cqfd

\vspace{3mm}

%  In order to achieve the proof of Theorem \ref{strongGWI} $(i)\Longrightarrow
%(ii)$, we still have 
It remains to consider the $x=0$ case: Let $t, s\in [0, \infty)$;
for any $p\geq 1$, we may and will choose  $s_p$ 
such that $|s_p -s| \leq 1/ \gamma_p$ and $[\gamma_p (s_p +t)]
=[\gamma_ps]+ [\gamma_p t]$. Use the Markov property for $Y^{*,p}$ at 
$[\gamma_p s]$ to get
\begin{equation}
\label{mmm}
 d_p(s_p +t,\lambda)=d_p(s, u_p(t,\lambda))+ b_p(t,\lambda). 
\end{equation}
Since
$b_p=d_p$, we get for any $t,s, \lambda, K \geq 0$, any $t\in T$ and any $p\geq 1$: 
\begin{equation}
\label{inegaal}
 d_p(s_p+t, \lambda ) - d_p(t, \lambda)\geq d_p(s, K \wedge u_p
(t,\lambda)). 
\end{equation}
%
% and then $b(\cdot ,\lambda)$ is non-decreasing and continuous. 
%Use monotonicity
%again to get for any 
%Then for any $t\in T$
%\begin{equation}
%\label{bmono}
% b_{p_k}(t, \cdot ) 
%\build{\longrightarrow}_{k\to\infty}^{\; } d(t, \cdot)
%\end{equation}
%uniformly on every compact subsets of the open interval 
%$(0, \infty)$. 
%Let $t\in T$ and
%$s\in [0, \infty)$. For any $p\geq 1$, we may and will choose  $s_p$ such that $|s_p -s| \leq 1/ \gamma_p$ and $[\gamma_p (s_p +t)]
%=[\gamma_ps]+ [\gamma_p t]$. Use the Markov property for $Y^{*,p}$ at 
%$[\gamma_p s]$ to get
%\begin{equation}
%\label{Markovimmdis}
% d_p(s_p +t,\lambda)=d_p(s, u_p(t,\lambda))+ b_p(t,\lambda).
%\end{equation}
%Since $d(\cdot ,\lambda)$ is continuous and since the $ b_{p_k}(\cdot , \lambda )$'s
%are non-decreasing, (\ref{Hellyly}) implies that $d(\cdot ,\lambda)$ is 
%non-decreasing. By a simple monotonicity argument, 
%$ d_{p_k}(s_{p_k} +t,\lambda)$ tends to $d(s+t,\lambda)$ when 
%$k$ goes to $\infty$. By (\ref{Hellyly}) and (\ref{egalili}), we get 
%$$ d(s+t,\lambda)=d(s,0)+ b(t,\lambda)= d(t,\lambda) \;, \quad t\in T \, , \,s\geq 0 \, , \, \lambda >0 .$$
%Since $T$ is dense and $d(\cdot , \lambda)$ is continous, we 
%deduce that $d(s,\lambda)=
%x\lambda $, $s \geq 0$ and $Z^*$ has 
%to be a constant process which contradicts 
%the assumptions of Theorem \ref{strongGWI} $(i)$. \cqfd
Observe next that since $b_p=d_p$, the $d_p(\cdot , \lambda)$'s and
$d(\cdot , \lambda)$ are non-decreasing. Since $d(\cdot , \lambda)$ is 
continuous, Dini's Theorem implies that $d_p(s_p +t,\lambda)$
tends to $d(s +t,\lambda)$ when $p$ goes to $\infty$. Choose $t$ in
$T$ 
and pass to the limit in (\ref{inegaal}) along the subsequence 
$(p_k , k\geq 0)$ to get 
$$ d(s+t, \lambda ) - d(t, \lambda)\geq d (s, K \wedge u (t,\lambda)),
\quad t\in T, \, s\geq 0, \, \lambda >0 .$$
%since
%$b_p=d_p$, then $d(\cdot , \lambda)$ is non-decreasing for any
%$\lambda>0$. 
%%
%
%Deduce from (\ref{Markovimmdis}) that for any $s, \lambda
%K \geq 0$, any $t\in T$ and any $p\geq 1$: 
%$$ d_p(s+t, \lambda ) - d_p(t, \lambda)\geq d_p(s, K \wedge u_p
%(t,\lambda)). 
%$$
%Then use (\ref{Helly3}) to get 
%$$ d(s+t, \lambda ) - d(t, \lambda)\geq d (s, K \wedge u (t,\lambda)). $$
If $u(t, \lambda ) =\infty $, then let $\lambda $  go to $0$ in the
previous inequality to get $d(s, K)=0$, for every $s,K\geq 0$, which
contradicts the assumption of Theorem \ref{strongGWI} $(i)$. Thus, 
$u(t, \lambda ) \neq \infty $ and
consequently ${\rm m}_t \neq 0$, for  
any $t\in T$. Then, use  (\ref{Helly3}) to pass to the limit in
(\ref{mmm}) and to get 
\begin{equation}
\label{marlim}
d(s+t, \lambda )=d(s, u(t,\lambda)) + d(t, \lambda) \; , 
\quad \lambda >0 , \, s\geq 0, \,  t\in T .
\end{equation}
Fix $s_0>0$. Since $\P (Z_{s_0} >0 )>0$, $d(s_0, \cdot)$ has to 
be a continuous increasing mapping. Then, it admits a 
continuous increasing inverse denoted by 
$\Delta$: $\Delta(d(s_0, \lambda))=\lambda$, $\lambda \geq 0$. We extend the definition
of $u$ on $[0, \infty )\times[0, \infty ) $ by setting
$$ u(t, \lambda):=\Delta \left( d(s_0+t, \lambda )-d(t, \lambda)\right) .$$
Then observe that $u(t, 0+)=0$ for any $t\in T$. 
Thus, ${\rm m}_t ([0, \infty ))=1$ , $t\in
T$. Now recall that for any fixed $\lambda\geq 0$
and any $p\geq 1$, $u_p (\cdot ,
\lambda)$ is monotone. It implies (\ref{Helly4}) by a standard monotonicity
argument. Now, use (\ref{cmono}) to deduce 
$\varphi( \lambda) <\infty$ , $\lambda >0$ and 
$\varphi( 0+)=0$, which both imply that ${\rm n}$ is a probability
measure on $[0, \infty)$. Then, use similar arguments to those used 
in the $x\neq 0$ case to complete the proof 
of Theorem \ref{strongGWI} $(ii)$. \cqfd

\vspace{3mm}

{\bf Proof of Theorem \ref{strongGWI} $(ii)\Longrightarrow
(iii)$}: The weak convergence of finite dimensional
marginals is a straightforward consequence of 
Lemma \ref{Grimvall} $(a) \Longrightarrow (b)$ combined with a simple computation 
based on (\ref{transGWI}). So, it remains to prove tightness. To that end, 
we adapt an argument of Grimvall \cite{Gr}: Fix $a\geq 0$ and 
denote by $(Y^{*,p}_n (a); n\geq 0)$ a GWI($\mu_p , \nu_p $)-process started 
at $Y^{*,p}_0 (a)=[pa]$.
%and take in particular $Y^{*,p}(x)=Y^{*,p}$. 
Denote by ${\bf Q}_a^{(p)} (\cdot )$ the distribution of $
p^{-1}( Y^{*,p}_{1} (a)-[ap])$. 
%$$ {\bf Q}_a^{(p)} (\cdot ) := \P \left( p^{-1}( Y^{*,p}_{1} (a)-[ap]) \,
%  \in \, \cdot \right) .$$
Theorem 2.2' \cite{Gr73} (see also Lemma 3.6 \cite{Gr}) asserts that
the sequence of the distributions of the processes 
$p^{-1}Y^{*, p}_{[\gamma_p \, \cdot \, ]} $, $p\geq 1$ is a tight
sequence in $\D ([0,
\infty ), \R )$ if the two following conditions are satisfied: 

\begin{description}
\item $(d)$ For any $t\geq 0$,  $ \quad \lim_{M\rightarrow \infty} 
\limsup_{p\rightarrow \infty} \P  \left( \sup_{0\leq s\leq t} 
 p^{-1}Y^{*,p}_{[\gamma_p s ]}>M \right) =0$.

\item $(e)$ For every compact set $C \subset [0,\infty )$, 
$$\{ ({\bf Q}_a^{(p)})^{*\gamma_p} \; , \; a\in C \, ,\,  p\geq 1 \} $$ 
is a tight family of probability measures on $\R$. 

\end{description}

\noindent
{\bf Proof of (e)}: Observe that 
$$ \left({\bf Q}_a^{(p)} \right)^{*\gamma_p}= \mu_p \left( \frac{\cdot
      -1}{p}\right) ^{*[ap]\gamma_p} *
\nu_p \left( \frac{\cdot
      }{p}\right) ^{*\gamma_p}  . $$
Thus, $(e)$ easily follows from Theorem \ref{strongGWI} $(ii)$. 
\cqfd

\vspace{3mm}

\noindent
{\bf Proof of (d)}: Fix $t>0$. Let $K$ be any positive real number. Observe
that for any
$p\geq 1$, any $\lambda, y >0$ and any $s\in [0,t]$, we have
\begin{eqnarray*}
\P \left(  p^{-1}Y^{*,p}_{[\gamma_p s ]} (y) \leq  K \right) & \leq & \exp \left(
  K -d_p(s, \lambda)\right) \\
& \leq &\exp \left(  K - \frac{ [py]}{p} u_p(s, \lambda) \right) \\
& \leq & \exp \left(  K - \frac{ [py]}{p} \lambda \wedge u_p(t, \lambda)
  \right)
\end{eqnarray*}
(use (\ref{1analy}) and the monotonicity of the 
$u_p (\cdot , \lambda)$'s). Now, since Theorem \ref{strongGWI} $(ii)$ implies Lemma
\ref{Grimvall} $(c)$, we get $\inf \{ \lambda \wedge u_p(t, \lambda) ,\; p\geq 1\} >0$. Thus, it proves
that for any $K>0$, there exists $M(K) >0$ such that 
\begin{equation}
\label{imiGrimest}
\P \left(  p^{-1}Y^{*,p}_{[\gamma_p s ]} (y) >  K \right) > 1/2 \; , \quad 
y\geq M(K)\,  , \; s\in [0,t] \, ,\; p\geq 1. 
\end{equation}
\noindent 
%{\bf Proof of Claim 4:} Let $(Y^{p}_n (y); n\geq 0)$
%stands for a GW($\mu_p$)-process 
%started at 
%$Y^{p}_0(y)=[py]$ and observe  
%that we can write $p^{-1}Y^{*, p}_{[\gamma_p s ]} (y)=A+B$ where $A$ and $B$ 
%are two non-negative independent random variables and where $A$ is 
%distributed as $p^{-1}Y^{p}_{[\gamma_p s ]} (y)$. It implies 
%\begin{equation}
%\label{imiGrimest2}
%\P \left(  p^{-1}Y^{*p}_{[\gamma_p s ]} (y) > K \right) \geq 
% \P \left(  p^{-1}Y^{p}_{[\gamma_p s ]} (y) > K \right) \; , \quad K >0 \,  ,
% \; s\in [0,t] \,  , \; p\geq 1.
%\end{equation}
%By Lemma \ref{Grimvall}, the first convergence of Theorem \ref{strongGWI}
%$(ii)$ implies that for any $y\geq 0$ 
%the following convergence  
%\begin{equation}
%\label{conauxi}
%\left( p^{-1}Y^{p}_{[\gamma_p t]} (y) \, ;\;  t\geq 0 \right) 
%\build{\longrightarrow}_{p\to\infty}^{{\rm (d)}} \left( Y_t (y) \,  ; \; t\geq 0 \right)
%\end{equation}
%holds weakly in 
%$\D( [0, \infty), \R)$ where $(Y_t(y); t\geq 0)$ stands for a CSBP($\psi$)
%started at $Y_0 (y)=y$. Since the 
%one-dimensional marginals of $Y(y)$ are stochastically increasing with $y$, 
%we can find $M(K)>0$ such that
%$$ \inf_{s\in [0,t]} \P \left( Y_s(y) > 2K  \right) >1/2 \; ,  \; \quad y\geq M(K) .$$
%Then, (\ref{imiGrimest}) follows from (\ref{imiGrimest2}) and  (\ref{conauxi}) \cqfd
%\vspace{3mm}

Now use the Markov property and (\ref{imiGrimest}) to get  
$$ \P \left(  p^{-1}Y^{*,p}_{[\gamma_p t ]} (x) > K  \right) \geq 
\frac{1}{2} \P \left( \sup_{0\leq s\leq t} 
 p^{-1}Y^{*,p}_{[\gamma_p s ]} (x) >M (K) \right) \; ,\quad p \geq 1$$
and $(d)$ follows from the one-dimensional marginals convergence. \cqfd

\subsection{Proof of Theorem \ref{convimmi}.}

We now consider a sequence of {\it (sub)critical} offspring distributions 
$(\mu_p;p\geq 1)$. Recall notations of Section 1 and 
consider a sequence $(\tau_p; p\geq 1)$
of GWI($\mu_p,r_p$)-trees where $\mu_p$ and $r_p$ satisfy (\ref{convXUV}) and
(\ref{extinctech}). 
For any $p\geq 1$, denote 
by $H^p=(H^p_k;k \geq 0)$ and by 
$D^p=(D^p_k;k\geq 0)$ resp. 
the height process and the random walk associated with the forest 
$f(\tau_p)$ containing the ``left bushes'' of $\tau_p$. 
We also denote by $H^{\bullet , p}$ and by 
$D^{\bullet , p}$ the processes corresponding to the forest 
$f(\tau_p^{\bullet})$ that contains the ``right bushes'' of $\tau_p$. 
Recall from Section 2.2 the notations $(L_k(\tau_p); k\geq 0)$
and $(L_k(\tau^{\bullet}_p); k\geq 0)$ , $p\geq 1$. 
We denote by $\Sigma $ the function such that  
 \begin{equation}
 \label{foncsig}
 \overleftarrow{H}^p = \Sigma (H^p, D^p, L(\tau_p)) \quad {\rm and } \quad  
 \overrightarrow{H}^p = \Sigma (H^{\bullet , p}, D^{\bullet , p}, 
 L(\tau^{\bullet}_p)) ,
 \end{equation}
that is specified by (\ref{spinaldec1}), 
(\ref{spinaldec2}), (\ref{spinaldec3}) and (\ref{spinaldec4}).

We use Corollary 2.5.1 \cite{DuLG} asserting that (\ref{convXUV}) and
(\ref{extinctech}) imply the joint convergence 
\begin{equation}
\label{jointconv}
\left(p^{-1}D^p_{[p\gamma_pt]},\gamma_p^{-1}H^p_{[p\gamma_pt]}
\,;\,t\geq 0\right)
\build{\longrightarrow}_{p\to\infty}^{{\rm (d)}} (X_t,H_t;\; t\geq 0)
\end{equation}
holds in distribution in $\D([0, \infty),\R^2)$. We also get an analogous 
convergence for $H^{\bullet , p}$ and $D^{\bullet , p}$ since they have the
same distribution.

%Next, observe that
%(\ref{convXUV}) and 
Next, we use Remark \ref{spinalGWI} and a standard argument 
to deduce from the right limit of (\ref{convXUV}) that the following
convergence 
\begin{equation} 
 \label{confoncLL}
 \left( p^{-1}L_{[\gamma_p t]}(\tau_p) \; , 
 \; p^{-1}L_{[\gamma_p t]}(\tau^{\bullet}_p) \, ; \, t\geq 0 \right)
 \build{\longrightarrow}_{p\to\infty}^{{\rm (d)}} (U_t, V_t 
 \, ; \, t\geq 0), 
 \end{equation}
holds in distribution in $\D([0, \infty), \R^2)$. Thus, 
by (\ref{jointconv}), (\ref{confoncLL}) 
and Skorohod's representation theorem, we may assume that the 
following convergences 
$$\lim_{p\rightarrow \infty }
 \left(p^{-1}D^p_{[p\gamma_pt]},\gamma_p^{-1}H^p_{[p\gamma_pt]}
 \right)_{t\geq 0}=(X,H) , $$ 
$$
\lim_{p\rightarrow \infty } 
 \left(p^{-1}D^{\bullet ,p}_{[p\gamma_pt]},\gamma_p^{-1}
 H^{\bullet ,p}_{[p\gamma_pt]}\right)_{t\geq 0} =(X',H') $$
and 
$$ \lim_{p\rightarrow \infty } \left( p^{-1}L_{[\gamma_p t]}(\tau_p) \; , 
 \; p^{-1}L_{[\gamma_p t]}(\tau^{\bullet}_p) \right)_{t\geq 0} = (U,V) $$
hold a.s. in $\D([0, \infty), \R^2)$, where $(X,H)$, $(X',H')$ and $(U,V)$ are
independent processes and where $(X,H)$ and $(X',H')$ have the same distribution.
For convenience of notation, we 
keep denoting in the same way the random processes involved in
the latter almost sure convergences, so we may also assume that (\ref{foncsig}) and (\ref{fonccontsig}) hold. We first prove 
\begin{equation}
\label{convleftH}
\left(\gamma_p^{-1}\overleftarrow{H}^p_{[p\gamma_pt]} ;\; t\geq 0
 \right)  
\build{\longrightarrow}_{p\to\infty}^{ \; }\left( \overleftarrow{H}_t ;
 t\geq 0 \right)
\end{equation}
a.s. in $\D([0, \infty), \R)$. To that end, let us introduce for any 
$\omega \in \D([0, \infty), \R)$, the right-continuous inverse of $\omega$ 
 $$ S_x (\omega )= \inf \{ s\geq 0 \; : \; \omega (s)>x \}  \; ,\quad x\in \R $$
(with the convention $\inf \emptyset =\infty$). Set $\v (\omega)= \{ x\in \R  \; : \; S_{x-} (\omega ) <
 S_x (\omega ) \}$. It is easy to check that $\omega \rightarrow S_x (\omega )$ is continuous in 
 $\D([0, \infty), \R)$ at any $\omega $ such that $x\notin  \v (\omega)$ 
 (see Proposition 2.11, Chapter VI \cite{Ja}). By (\ref{UVcont}), the process
 $x\rightarrow S_x (U)=U^{-1}_x$ has a.s. continuous
 sample paths. Then it implies a.s. 
 $$\lim_{p\rightarrow \infty }S_x \left( p^{-1}
 L_{[\gamma_p \, \cdot ]} (\tau_p) \right)= U^{-1}_x  \; , \; \quad x\in \Q_+ .$$
 Since $U^{-1}$ is a continuous increasing process, standard arguments imply  
 \begin{equation}
 \label{convinvU}
 \left ( S_x ( p^{-1} L_{[\gamma_p \, \cdot ]}(\tau_p) ) 
\, ; \; x\geq 0 \right)
 \build{\longrightarrow}_{p\to\infty}^{ \; } \left( U^{-1}_x ; x\geq 0\right) 
 \end{equation}
 a.s. in $\D([0, \infty), \R)$ (see
 Theorem 2.15, Chapter VI \cite{Ja}). Let us set
 $$ \aal^p (n):=\inf \{ k\geq 0 :\; L_k(\tau_p)  \geq 1- \inf_{j\leq n } D^p_j 
 \}  \; , \; \quad p,n\geq 1.$$
 Since $t\rightarrow I_t=\inf_{s\in [0,t]} X_s$ is a continuous process, the following 
 convergence 
 $$ \left( \inf_{j\leq [p\gamma_pt]} p^{-1} D^p_j \, ; \;  t\geq 0\right) 
\build{\longrightarrow}_{p\to\infty}^{{\rm (d)}} \left(I_t ; t\geq 0\right) $$
 a.s. holds uniformly on every compact subsets of $[0, \infty )$. Thus, by (\ref{convinvU}) 
 \begin{equation}
 \label{convalpha}
 \left( \gamma_p^{-1}\aal^p ( [p\gamma_pt]) ; t\geq 0\right) \build{\longrightarrow}_{p\to\infty}^{{\rm (d)}}  
 (U^{-1}_{-I_t} ;\; t\geq 0)
 \end{equation}
a.s. uniformly on every compact subsets of $[0, \infty )$. Next, we set  
$$ {\bf \pi}^p (n) := \inf \{ k\geq 0\; :\; k+\aal^p (k) \geq n \} \; ,\;
 \quad p,n \geq 1.$$
 Then, by (\ref{spinaldec3}) and 
 (\ref{spinaldec4}), we get 
 \begin{equation}
 \label{identite}
 \overleftarrow{H}^p_n =n- {\bf \pi}^p (n) + H^p_{{\bf \pi}(n)} \quad 
 {\rm and } \quad \aal^p ( {\bf \pi}^p (n) -1) \leq n- {\bf \pi}^p (n) 
 \leq \aal^p ({\bf \pi}^p (n)) .
 \end{equation}
 We easily deduce from (\ref{convalpha}) that 
 $(p\gamma_p)^{-1}{\bf \pi}^p ([p\gamma_p \cdot])$  a.s. converges to the identity
 map uniformly on every compact subsets of $[0,\infty)$. 
Then, (\ref{convleftH}) follows from (\ref{identite}) and (\ref{convalpha}). Use
 similar arguments to prove the corresponding convergence for $\overrightarrow{H}^p$. 
 Thus, we have proved that the following convergence 
 \begin{equation}
 \label{functionconv}
 \left( \gamma_p^{-1} \overleftarrow{H}^p_{[p\gamma_pt]}, 
 \gamma_p^{-1} \overrightarrow{H}^p_{[p\gamma_pt]} \,;\,t\geq 0\right)
 \build{\longrightarrow}_{p\to\infty}^{{\rm (d)}} (
 \overleftarrow{H}_t , \overrightarrow{H}_t ;\; t\geq 0)
 \end{equation}
 holds in distribution in $\D ([0, \infty) , \R^2)$. The joint convergence of the
 corresponding contour processes is a consequence of (\ref{contour1}) and
 (\ref{contour2}): denote by $q_p$ the function 
 associated with $\tau_p$ as defined in Section 2.2; set 
 $y_p(s)=(p\gamma_p)^{-1}
 q_p(p\gamma_p s)$; by (\ref{contour1}) 
%and (\ref{functionconv}) 
we get for every $T>0$,
 \begin{equation}
 \label{cont-tech1}
 \sup_{s\leq T}\Big|{1\over \gamma_p}
 \overrightarrow{C}_{p\gamma_p s}(\tau_p) -{1\over \gamma_p}
 \overrightarrow{H}^p_{p\gamma_p y_p(s)}\Big|
 \leq{1\over \gamma_p}+{1\over \gamma_p}\sup_{n\leq T p\gamma_p}
 |\overrightarrow{H}^p_{n+1}-\overrightarrow{H}^p_{n}|
 \build{\longrightarrow}_{p\to\infty}^{} 0 
 \end{equation}
 in probability by (\ref{functionconv}). On the other hand, we get from
 (\ref{contour2}) 
%and (\ref{functionconv})
 \begin{equation}
 \label{cont-tech3}
 \sup_{s\leq T}|y_p(s)-{s\over 2}|
 \leq{1\over 2p\gamma_p}\sup_{k\leq T p\gamma_p}\overrightarrow{H}^p_k
 +{1\over p\gamma_p}
 \build{\longrightarrow}_{p\to\infty}^{} 0 
 \end{equation}
 in probability by (\ref{functionconv}). Using similar arguments, we 
 prove analogue convergences in probability for the right 
 contour processes. Then, 
\begin{equation}
 \label{functionconvcont}
 \left( \gamma_p^{-1} \overleftarrow{C}_{2p\gamma_pt}(\tau_p), 
 \gamma_p^{-1} \overrightarrow{C}_{2p\gamma_pt}(\tau_p ) \,;\,t\geq 0\right)
 \build{\longrightarrow}_{p\to\infty}^{{\rm (d)}} (
 \overleftarrow{H}_t , \overrightarrow{H}_t ;\; t\geq 0).
\end{equation}

  The convergence of the sequence 
$p^{-1} Y^{*,p}_{[p\gamma_p \, \cdot \,]} $ is a consequence of 
Theorem \ref{strongGWI}: we easily see that  (\ref{convXUV})
implies (\ref{hypotheses}) by taking $W_1= U_1+V_1$ 
and thus,  
%$\varphi (\lambda )= \Phi (\lambda , \lambda)$) and by taking 
$$ \nu_p (k-1)= \sum_{1 \leq j\leq k} r_p(k,j) \; ,\quad k\geq 1 .$$
To get the desired joint convergence of the contour processes with the 
GWI process, argue exactly as in the proof
of corollary 2.5.1 \cite{DuLG} p. 63-64. \cqfd

%\begin{center}

%\includegraphics{GWI_1.pdf}

%\includegraphics[ viewport=10 150 160 220]{GWI_1.pdf}

%\end{center}

%\begin{center}

%   \includegraphics{sin_1.pdf}
  
%   \includegraphics[viewport=10 80 160 220]{sin_1.pdf}

%\end{center}

\bibliographystyle{plain}

\end{document}